\documentclass{article}
\usepackage{theorem}    % for the \theoremstyle, cf. Samarin
\usepackage{amsmath}
\usepackage{amsfonts}
\usepackage{version}
\usepackage{color}      % xfig pishet \color{}
\usepackage{graphicx}   % ps pictures
\usepackage[all]{xy}
\numberwithin{equation}{section}
\theorembodyfont{\rmfamily}

\newtheorem{thing}{}[subsubsection]   %!!!!! 
\theoremstyle{plain}

\newtheorem{proposition}{Proposition}
\newcommand{\al}{\alpha}
\newcommand{\be}{\beta}
\newcommand{\fibered}{\mathop{\times}\limits}  % fibered product 
\newcommand{\dual}{^{\vee}}

\newcommand{\jk}{{J_{\xi}}}
\newcommand{\ok}{{\OO_{\xi}}}
\newcommand{\ox}{{\OO_{X}}}
\newcommand{\oc}{{\OO_{C}}}

\newcommand{\OO}{{\cal O}} % Structure sheaf
\newcommand{\PP}{{\mathbb P}} % projective space
\newcommand{\Spec}{\mathop{\rm Spec}}
\newcommand{\Aut}{\mathop{\rm Aut}}
\newcommand{\Tors}{\mathop{\rm Tors}} % torsion

\newcommand{\Num}{\mathop{\rm Num}}

\newcommand{\Pic}{\mathop{\rm Pic}\nolimits}

\newcommand{\triv}[1]{#1 \tensor \OO_X}   % V tensor \OO

\newcommand{\ZZ}{\mathbb Z}
\newcommand{\QQ}{\mathbb Q}

\newcommand{\rk}{\mathop{\rm rk}}

\newcommand{\codim}{\mathop{\rm codim}\nolimits}
\newcommand{\Deg}{\mathop{\rm Deg}}

\newcommand{\Ext}{\mathop{\rm Ext}\nolimits}  
\newcommand{\Hom}{\mathop{\rm Hom}\nolimits}  
\newcommand{\Monom}{\mathop{\rm Monom}\nolimits}  
\newcommand{\RHom}{\mathop{\rm RHom}\nolimits}
\newcommand{\Tor}{\mathop{\rm Tor}\nolimits}  
\newcommand{\sheafExt}{\mathcal E \mathit{xt}}       %curved Ext
\newcommand{\sheafHom}{\mathcal H \mathit{om}}
\newcommand{\coker}{\mathop{\rm coker}}
\renewcommand{\Im}{\mathop{\rm Im}}

\newcommand{\End}{\mathop{\rm End}}
\newcommand{\depth}{\mathop{\rm depth}}

\newcommand{\tensor}{\mathop{\otimes}}

\newcommand{\Supp}{\mathop{\rm Supp}}  
  
\newcommand{\ch}{\mathop{\rm ch}}  
\newcommand{\Gr}{\mathop{\rm Gr}}      
\newcommand{\id}{{\rm id}}

\newcommand{\intersect}{\cap}
\newcommand{\intersectIn}[1]{\mathop{\cap}\limits_{#1}}
\newcommand{\union}{\bigcap}
\newcommand{\isom}{\stackrel{\sim}{\to}}

\newcommand{\strelka}[1]{\stackrel{#1}{\to}} % arrow s nadpis'yu

\newcommand{\triple}[3]{
                         0 \to {#1} \to {#2} \to {#3} \to 0 
                       }

\newcommand{\twoVector}[2]{%
  \begin{bmatrix}
     {#1} \\ {#2}
  \end{bmatrix}
}

\newcommand{\xyArrow}[2]{\ar@{}[#2]|-*{#1}}
\newcommand{\xyIsom}{\ar^-*@{~}}          % isomorfism: strelka s ~ naverhu,
\newcommand{\xyIsomBelow}{\ar_-*@{~}}     % to zhe, no ~ vnizu
\newcommand{\xySimeq}{\ar@{-}^*@{~}}      % \simeq - strelka 
\newcommand{\xyRavno}{\ar@{=}}
\newcommand{\xyMaps}{\ar@{|->}}           % element perehodit v element
\newcommand{\xyKrivayastrelka}{\ar@{~>}}  % krivaya strelka
\newcommand{\xyMono}{ \ar@{^{(}->}}       %vlozhenie
\newcommand{\xyMonom}{\ar@{_{(}->}}       %tozhe vlozhenie
\newcommand{\xyEpi}{\ar@{->>}}            %epimorfizm
\newcommand{\xyCherta}{\ar@{-}}           %pryamaya liniya
\newcommand{\xySubsetSmall}{\ar@{}|(0.60){\subset}} % malen'kaya strelka subset
\newcommand{\xySubset}{\ar@{}|(0.60)*{\subset}} % subset
\newcommand{\xySubsetA}{\ar@{}|*{\subset}[r]} % subset, raspolozhennyj
\newcommand{\xyCap}{\ar@{}|(0.60)*{\cap}} % cap (vlozhenie vniz)
\newcommand{\xyCapA}{\ar@{}|*{\cap}[d]}   %      vlozhenie vniz, raspolozheno
\newcommand{\xyCup}{\ar@{}|*{\cup}[u]}    %      vlozhenie vverh 

\newcommand{\smallXyMatrix}{\xymatrix@R=10pt@C=10pt}

\newcommand{\xyTriple}[5]{%
   \smallXyMatrix
   {
      0      \ar[r]        & 
      {{#1}} \ar[r]^-{{#4}} & 
      {{#2}} \ar[r]^-{{#5}} & 
      {{#3}} \ar[r]        & 
      0 
   }
}

\newcommand{\displayD}[7]{
  \smallXyMatrix
  {
    {} & 0\ar[d] & 0\ar[d] & {} & {} \\
    {} & {#1}\xyRavno[r]\ar[d] & {#1}\ar[d] & {} & {} \\
    0\ar[r] & {#2}\ar[r]\ar[d] & {#3}\ar[r]\ar[d] & {#4}\ar[r]\xySimeq[d] & 0\\
    0\ar[r] & {#5}\ar[r]\ar[d] & {#6}\ar[r]\ar[d] & {#7}\ar[r] & 0\\
    {} & 0 & 0 & {} & {} 
  } 
}

\newcommand{\displayF}[7]{
  \smallXyMatrix
  {
    {} & 0\ar[d] & 0\ar[d] & {} & {} \\
    {} & {#1}\xyRavno[r]\ar[d] & {#1}\ar[d] & {} & {} \\
    0\ar[r] & {#2}\ar[r]\ar[d] & {#3}\ar[r]\ar[d] & {#4}\ar[r]\xyIsom[d] & 0 \\
    0\ar[r] & {#5}\ar[r]\ar[d] & {#6}\ar[r]\ar[d] & {#7}\ar[r] & 0 \\
    {} & 0 & 0 & {} & {} 
  } 
}

\newcommand{\displayE}[6]{
  \xymatrix
  {
    {} & 0\ar[d] & 0\ar[d] & {} & {} \\
    {} & {#1}\xyRavno[r]\ar[d] & {#1}\ar[d] & {} & {} \\
    0\ar[r] & {#2}\ar[r]\ar[d] & {#3}\ar[r]\ar[d] & {#4}\ar[r]\xyRavno[d] & 0\\
    0\ar[r] & {#5}\ar[r]\ar[d] & {#6}\ar[r]\ar[d] & {#4}\ar[r] & 0\\
    {} & 0 & 0 & {} & {} 
  } 
}
\newcommand{\displayC}[7]
{
  \smallXyMatrix
  {
     {} & {} & 0 \ar[d] & 0 \ar[d] & {} \\
     {} & {} & {{#1}} \xyRavno[r] \ar[d]^{{#7}} & {{#1}} \ar[d] & {} \\
     0 \ar[r] & {{#2}} 
                \xyRavno[d] 
                \ar[r]      & {{#3}}
                              \ar[d]\ar[r] & {{#4}}\ar[d]\ar[r] & 0 \\
   0 \ar[r] & {{#2}} \ar[r] & {{#5}} \ar[d]\ar[r] & {{#6}} \ar[d]\ar[r] & 0\\
   {} & {} & 0 & 0 & {} 
  }
}

\newcommand{\displayCa}[9]
{
  \smallXyMatrix
  {  
     {} & {} & 0 \ar[d] & 0 \ar[d] & {} \\
     {} & {} & {{#1}} \xyIsom[r] \ar[d]^{{#8}} & {{#2}} \ar[d]^{{#9}} & {} \\
     0 \ar[r] & {{#3}} 
                \xyRavno[d] 
                \ar[r]      & {{#4}}
                              \ar[d]\ar[r] & {{#5}}\ar[d]\ar[r] & 0 \\
   0 \ar[r] & {{#3}} \ar[r] & {{#6}} \ar[d]\ar[r] & {{#7}} \ar[d]\ar[r] & 0\\
   {} & {} & 0 & 0 & {} 
  }
}
\newcommand{\correspondence}[5]
{ 
   \smallXyMatrix
   {
        {} & {{#1}}\ar[ld]_-{{#4}}\ar[rd]^-{{#5}} & {} \\
        {{#2}} & {} & {{#3}}
   }
}   
    
\newcommand{\xyEquivalenceA}[5]
{
   \smallXyMatrix
   {
     0 \ar[r] & {#1} \ar[r]\xyRavno[d] & {#2} \ar[r]\xyIsomBelow[d]^{#5} & {#3} \ar[r]\xyRavno[d] & 0 \\
     0 \ar[r] & {#1} \ar[r] & {#4} \ar[r] & {#3} \ar[r] & 0
   }
}

\newcommand{\xyInducedOnTheRight}[7]
{
   \smallXyMatrix
   {
     0\ar[r] & {#1}\ar[r]\xyRavno[d] & {#2}\ar[r]\ar[d]^{#6} & {#3}\ar[r]\ar[d]^{#7} & 0 \\
     0\ar[r] & {#1}\ar[r] & {#4} \ar[r] & {#5} \ar[r] & 0
   }
}
\newcommand{\xyInducedOnTheLeft}[7]
{
   \smallXyMatrix
   {
     0 \ar[r] & {#1} \ar[r]\ar[d]^{#6} & {#2} \ar[r]\ar[d]^{#7} & {#3} \ar[r]\xyRavno[d] & 0 \\
     0 \ar[r] & {#4} \ar[r] & {#5} \ar[r] & {#3} \ar[r] & 0 \\
   }
}
\excludeversion{for-me}

\begin{document}

\title{On the Brill-Noether theory for K3 surfaces II }
\author{Maxim Leyenson}
% \date{}

\maketitle
%-------------------------------------------------
\begin{abstract}
    Let $(S,H)$ be a polarized K3 surface, $E$ be a coherent sheaf on
$S$ and $W \subset H^0(S,E)$ be a linear subspace. If we are lucky,
there is an exact sequence
$$0 \to W \tensor \OO_S \to E \to E' \to 0 $$ which gives a
correspondence between moduli spaces of sheaves of different ranks
which we used in the first part of the paper to establish some
properties of Brill-Noether loci in the moduli space. We allow $E$ to
be locally free, torsion free of rank one or a line bundle with
support on a curve, thus studying simultaneously Brill-Noether special
vector bundles, special 0-cycles and special linear systems on curves.

   To complete the work begun in the Part 1, we need to establish a
number of properties of this correspondence. In this paper we prove
that it behaves nicely for globally generated vector bundles,
establish the existence of globally generated vector bundles in moduli
spaces on K3, and prove that the correspondence preserves stability if
$\Pic S = \ZZ c_1(E)$, thus completing the work of Part 1.
\end{abstract}
%-------------------------------------------------

\tableofcontents

\subsection{Introduction} 
Let $X$ be an algebraic surface. (Later we will need to assume that
the characteristics of the base field is 0.)

We call a coherent sheaf $E$ on $X$ acceptable if $\rk E \ge 2$ and
$E$ is locally free, or if $\rk E = 1$ and $E$ is torsion free, or if
$\rk E = 0$ and $E$ is an invertible sheaf on a curve on $X$.

   Let $A$ be the moduli space of pairs $(E,V)$, where $E$ is a stable
acceptable sheaf of a given class in $K_0(S)$ and $V \subset H^0(X,E)$
(it is constructed later in the paper). If we are lucky, a point in
this moduli space gives an exact sequence of the form
$$0 \to V \tensor \OO_X \to E \to E' \to 0,$$ where $E'$ is also
acceptable and stable. This gives a correspondence between moduli
spaces of acceptable sheaves of different ranks which we studied in
~\cite{part1}.

It turns out to be non-trivial to prove that this correspondence is
non-empty and well-behaved. In this paper we prove that it behaves in
``naively expected way'' for globally generated vector bundles in the
case $\Pic X = \ZZ$. In particular, we prove the following results we
use in ~\cite{part1}:

\begin{thing} {\bf Monomorphic evaluation maps:}
 Let $X$ be an algebraic variety, $E$ be a globally generated vector
bundle on $X$, and $l \le \rk E$. Then for a {\it generic} vector
space $V \in Gr(l,H^0(X,E))$ the evaluation map $V \tensor \ox \to E$
is monomorphic.
\end{thing}

\begin{thing} {\bf Acceptable extensions:} 
   Let $X$ be an algebraic surface, $E'$ be an acceptable sheaf on $X$,
and let $W$ be a $k$-vector space. Let $\mathbb A = \Ext^1_X(E',W
\tensor \OO_X)$. A point $e \in \mathbb A$ gives an extension class of
the form
  \begin{equation*}
      \triple{W \tensor \OO_X}{E}{E'}
   \end{equation*}
   \begin{enumerate}
      \item If $\rk E' \ge 2$, then $E$ is acceptable (this is a
   trivial part);
      \item if $\rk E' = 1$, $E' = \jk(L)$, where $L$ is ample,
   and the pair $(\xi,L+K_X)$ is Caley-Bacharash, then for a {\it generic} 
   $e \in \mathbb A$ the sheaf $E$ is acceptable;
      \item if $\rk E' = 0$, $E'=(i_C)_* (B)$, and $A :=
   N_{C/X}B^{-1}$ is globally generated, then for a {\it generic} $e
   \in \mathbb A$ the sheaf $E$ is acceptable.
   \end{enumerate}
\end{thing}

\begin{thing} {Acceptable factors:}
   Let $X$ be a smooth algebraic surface in characteristics 0, $E$ be
   a globally generated vector bundle on $X$, and $1 \le l \le \rk E$.
   Then for a {\it generic} vector space $V \in Gr(l,H^0(E))$ the
   cokernel of the evaluation map $e_V: V \tensor \OO_X \to E$ is
   acceptable.
\end{thing}

\begin{thing}{\bf Stable extensions:}
   Let $(X,H)$ be a polarized algebraic surface, and let $E'$ be an
   $H$-stable acceptable sheaf on $X$. Assume that we are given an
   extension
   \begin{equation*}
      0 \to W \tensor \OO_X \to E \to E' \to 0
   \end{equation*}
   Let $e$ be the class of this extension in $\Ext^1(E',W \tensor \OO_X)$,
   and let $\alpha_e$ be the image of $e$ under the isomorphism
   $\Ext^1(E',\triv{W}) \isom \Hom(W\dual,\Ext^1(E',\OO))$.

   \begin{enumerate}
     \item If $\alpha_e$ is not injective, then $E$ is not $H$-stable;
     \item If $\alpha_e$ is injective, $E'$ is acceptable, and 
           $\Pic X \simeq \ZZ \cdot c_1(E)$, where $c_1(E)$ is ample, 
           then $E$ is $H$ - stable.
   \end{enumerate}
\end{thing}

\begin{thing}{\bf Stable factors:}
  Let $(X,H)$ be a polarized algebraic surface, $E$ be a stable
coherent sheaf on $X$, and assume that we are given an exact sequence
of the form
\begin{equation*}
   0 \to W \tensor \OO_X \to E \to E' \to 0
\end{equation*}
where $E'$ is torsion-free.  Assume that $\Pic X \simeq \ZZ \cdot
c_1(E)$, where $c_1(E)$ is ample. Then $E'$ is $H$-stable.
\end{thing}

\begin{thing}{\bf Existence of globally generated vector bundles in
moduli spaces on K3.} Let $X$ be a K3 surface (this is the only place
in this paper where we use that $X$ K3). Assume that $\Pic X \simeq
\ZZ h$, $h$ ample, and consider a nonempty moduli space $M$ of vector
bundles on $S$ with $\rk E = r$, $c_1(E) = h$ and $c_2(E) = d$. Then
there is a globally generated vector bundle $E \in M$.
\end{thing}

  We also construct explicit complex of locally free sheaves computing
the cohomology groups of a locally free sheaf varying in a family and
prove numerous results on special 0-cycles, linear systems and vector
bundles.

\subsection{Acknowledgments.}
In addition to all the acknowledgements expressed in the previous
paper I would like to thank Andrey Levin (ITEP, Moscow) for very
helpful discussions. The idea of the proof of Lemma
~\ref{grassmanian-lemma} given in this paper belongs to him.

% -----------------------------------------------------------------

\section{Preliminaries.}

\input{preliminaries.duality}             
           %no [], spell-checked
           %\subsection{(.) Two duality lemmas.}

\subsection{ Acceptable sheaves on surfaces}  
{\bf Definition.} 
Let $X$ be an algebraic surface.  We say that a coherent sheaf $E$ is {\it
acceptable} if\\

\begin{tabular}{ll}
   $\rk E \ge 2$  & and $E$ is locally free; \\
   $\rk E = 1$    & and $E$ is torsion free; \\
   $\rk E = 0$    & and $E$ is a direct image of an invertible sheaf on
a Cartier divisor on $X$.
\end{tabular}

% perhaps for rk = 1 we want to add ``xi is a locally complete intersection''?

\subsection{ One corollary of the Hilbert-Burch theorem}

{\bf Lemma.} Let $R$ be a CM ring of dimension two and let $\phi: R^m
\to R^{m+1}$ be a monomorphism of free modules. If $\phi \tensor k(p)$
is monomorphic outside a codimension two subscheme in $\Spec R$, then
$\coker \phi$ is isomorphic to an ideal $J \subset R$ of depth 2.

\begin{for-me}
{\bf Proof.} Let $\xi$ be the degeneration scheme of $\phi$. Then
$\codim \xi = 2$. Since $X$ is CM, we have $\depth J_{\xi}$ = 2 and
the Hilbert-Burch theorem is applicable.
\end{for-me}

\subsection{ Poincare polynomials}

Let $X$ be a scheme and $E$ be a coherent sheaf on $X$. Let $\ox(1)$
be an ample invertible sheaf on $X$. The Poincare polynomial of $E$ is
\begin{multline*}
    P(E) = \chi(E(m)) = 
            \chi(E) + \chi(E|_H) m + \chi(E|_{H^2}) \frac{m(m+1)}{2} +
            \dots \\
                   + \chi(E|_{H^d}) \frac{m(m+1)\dots(m+d-1)}{d!},
\end{multline*}
where $H \in |\ox(1)|$ is a generic divisor in the linear system
$\ox(1)$, $H^i$ is the intersection of $i$ generic divisors in the
linear system $\ox(1)$, and $d$ is the dimension of support of $E$ 
(~\cite{huybrechts-lehn}).

In particular, if $\dim X = 2$, 
\begin{gather*}
       P(E) = \chi(E|_{H^2}) \frac{m(m+1)}{2} + \chi(E|_H) m + \chi(E) =   \\
              r (H,H) \frac{m(m+1)}{2} + \left( r(1-g) + (c_1(E),H) \right) m 
              + \chi(E)
\end{gather*}
where $r = \rk E$ and $g$ is genus of the curve $H$.

If $r > 0$, then the reduced Poincare polynomial is defined as 
$$
       p(E) = P(E)/r =
        (H,H) \frac{m(m+1)}{2} 
      + \left( (1-g) + \frac{(c_1(E),H)}{r} \right) m
      + \frac{\chi(E)}{r} 
$$

\subsection{ Vanishing of $H^2(X,E)$ for stable $E$ with positive $c_1$ }

{\bf Lemma.} Let $(X,\ox(1))$ be a polarized K3 surface and
$E$ be a stable coherent sheaf on $X$ of rank $r > 0$ such that
$(c_1(E),H) > 0$. Then $H^2(X,E) = 0$.

{\bf Proof.} By duality, $H^2(X,E) \simeq \Ext^0(E,\Omega^2_X)\dual =
\Hom(E,\ox)\dual$.

Assume that $\Hom(E,\ox) \ne 0$, i.e., that there is a nonzero homomorphism
$\phi: E \to \ox$. Let $J = \phi(E)$; $J$ is a sheaf of ideals of some closed
subscheme $Z \subset X$. The map $\phi$ induces an epimorphism
$$
   E \to J \to 0
$$
and therefore the stability of $E$ should imply $p(E) < p(J)$, i.e.,
$$
   \frac{(c_1(E),H)}{r} m + \frac{\chi(E)}{r} <
   c_1(J|_H) m + \chi(J)
$$
% (If $H$ is generic enough), 
One can see from the exact sequence
$$
   \triple{J \tensor \OO_H}{\OO_H}{\OO_{Z \intersect H}}
$$
that $c_1(J|_H) = - c_1(\OO_{Z \intersect H}) < 0$, 
while $(c_1(E),H) > 0$, which gives a contradiction.

%           no [], spell-checked
%           \subsection{Acceptable sheaves on a surface} proverit' 
%           \subsection{(.) One corollory of the Hilbert-Burch theorem}
%           \subsection{(.) Poincare polynomials}
%           \subsection{(.) Vanishing of $H^2(X,E)$ for a stable $E$ on K3}

\input{ structure-of-extensions.part1  }   
          % no [], spell-checked
          %\subsection{ (.) Extensions of modules over a commutative ring.}
          % \subsection{ (.) $\sheafExt^1(\ok,\OO_X) = 0$ } -
          % commented out
          % \subsection{ (.) Computation of $\sheafExt^i(\jk(L),\OO_X)$ }
          % \subsection{(.) Computation of $H^0(X,\sheafExt^1(\jk(L),\OO_X))$ }
          % \subsection{(.)  Spectral sequence for $\Ext^1(F,G)$}
          % \subsection{(.) Spectral sequence for $\Ext^1(F,\OO)$ 

\providecommand{\al}{\alpha}
\providecommand{\be}{\beta}
\providecommand{\ab}{\alpha\beta}
\providecommand{\bc}{\beta\gamma}
\providecommand{\ac}{\alpha\gamma}
\providecommand{\Ua}{U_{\alpha}}
\providecommand{\Aa}{A_{\alpha}}
\providecommand{\Fa}{F_{\alpha}}
\providecommand{\Ea}{E_{\alpha}}
\providecommand{\Ga}{G_{\alpha}}
\providecommand{\ga}{g_{\alpha}}
\providecommand{\ea}{e_{\alpha}}
\providecommand{\Eb}{E_{\beta}}
\providecommand{\gb}{g_{\beta}}
\providecommand{\phiab}{\phi_{\ab}}
\providecommand{\psiab}{\psi_{\ab}}
\providecommand{\Aab}{A_{\ab}}
\providecommand{\Uab}{U_{\ab}}
\providecommand{\Fab}{F_{\ab}}
\providecommand{\tg}{\tilde g}

\subsection{ Locally free extensions $\triple{\OO}{E}{\jk(L)}$}
\subsubsection{ Local case}
% Extensions $\triple{A}{E}{m}$ over a regular local ring of dimension two
\label{local-freeness-over-a-local-ring}
{\bf Lemma.} Let $A$ be a regular local ring of dimension two,
$m$ be its maximal ideal, and let 
\begin{equation}
   \triple{A}{E}{m}     
   \tag{$*$}
\end{equation}
be some extension of $A$-modules with class 
$e \in \Ext^1_A{m,A}$.
Then $E$ is a free module iff $e \ne 0$.

This is a classical result, for another proof see, f.e.,
~\cite{griffiths-harris}, p.5.4.

{\bf Proof.}
 If $e = 0$, then $E$ is not free, so assume $e \ne 0$. 
Let $p$ be $m$ considered as a point in $\Spec A$.
If $\{x,y\}$ is a regular sequence generating $m$, then $m$ has a Koszul
resolution %$$ \xyTriple{A}{A^2}{m}{ \twoVector{x}{y} }{ [-y,x] }
$$ \triple{M_1}{M_0}{m} $$
where $M_1 = A$ and $M_o = A^2$ 
and the exact triple (*) has a resolution
\begin{equation}
   \displayCa{E_1}{M_1}{A}{E_0}{M_0}{E}{m}{v}{w}
   \label{diagram:resolution}
\end{equation}
where $E_0 = E \mathop{\oplus}\limits_m M_0 = \ker(E \oplus M_0 \to m)$
is a free module of rank three.
% (it follows from $\Tor_1(E_0,k)=0$)

%The long exact sequence
%$$
%   0 \to \Tor_1(E,k(p))
%     \to \Tor_1(m,k(p))
%     \stackrel{\delta_1}{\to}
%     k(p)
%$$
%associated with (*) 
%shows that $E$ is free iff $\delta_1$ is injective. 
%
% mozhno prosche!

The middle column of ~\ref{diagram:resolution} shows that $E$ is free
iff $v(p) \stackrel{def}{=} v \tensor k(p) \ne 0$.

Now ~\ref{diagram:resolution}
generates
$$
  \smallXyMatrix
  {
     {} & {} & {0} & {0} & {} \\
     {} & {0} & {\Hom(E_1,A)} \ar[l]\ar[u] & {\Hom(M_1,A)}\xyIsom[l]\ar[u] & {0} \ar[l] \\
     {0} & {\Hom(A,A)}\ar[l]\ar[u] & {\Hom(E_0,A)}\ar[l]\ar[u]^{v^*} & {\Hom(M_0,A)}\ar[l]\ar[u]^{w^*} & {0} \ar[l]
  }
$$

Let $\id_A \in \Hom(A,A)$ be the identity homomorphism and let
us choose some lift $\tilde{\id}_A \in \Hom(E_0,A)$  of $\id_A$
to $\Hom(E_0,A)$. Then, by definition, 
$$  
    e = a^* ( v^* (\tilde{\id}_A) ) \mod \Im w^*
$$
Let $\phi = v^* (\tilde{\id}_A)$.
There is an isomorphism $\Hom(M_1,A) \simeq A$ which takes
$\Im w^*$ to $m$, so $    e = a^* (\phi) \mod m $.

  But $\tilde{\id}_A$ is the left splitting of the exact sequence
in the middle row of ~\ref{diagram:resolution}, 
as in the
$$
   \xymatrix
   {
     {} & {} & {E_1} \xyIsom[r]\ar[d]^v \ar@{..>}[ld]_{\phi} & {M_1}\ar[d] & {}\\ 
     {0}\ar[r] & {A}\ar[r]_{i} & {E_0} \ar[r]\ar@<-1mm>@{..>}[l]_{\tilde{\id}} & {M_0}\ar[r] & {0}
   }
$$

Since $v=i \circ \phi$  $v(p)=i(p) \circ \phi(p)$. But $i(p) \ne 0$
and $e\ne 0$ implies $\phi(p) \ne 0$ , which implies $v(p) \ne 0$.

\subsubsection{  Global case }
% Local freeness of extensions $\triple{\OO}{E}{\jk(L)}$
\label{local-freeness}
We will call an extension 
$$ \triple{\OO}{E}{\jk(L)} $$ locally free
iff $E$ is locally free.

{\bf Lemma.}
  Let $X$ be a smooth surface, $L$ be an ample line
bundle on $X$, and $\xi$ be a simple effective 0-subscheme in $X$. Then

   {\bf (a)} $\dim \Ext^1(\jk(L),\OO) = h^1(\jk(LK)) = \delta(\xi,LK)$,

   {\bf (b)} Let $\dim_k \Ext^1(\jk(L),\OO) > 0$ and $p \in \xi$.
Then either

   every curve in the linear series $|L+K|$ passing through $\xi - p$
also contains $p$, and the set of extensions parametrized by
$\Ext^1(\jk(L),\OO)$ and locally free at $p$ is a complement to a
hyperplane in $\Ext^1(\jk(L),\OO)$, 

or

there is a curve in the linear system $|L+K|$ containing $\xi-p$ but
not $p$, and all extensions classes in $\Ext^1(\jk(L),\OO)$ are not
locally free at $p$.

\begin{for-me}
~\ref{generic-extension- of-a-generic-special-J-is-locally-free}
pokazyvaet, chto pervoe uslovie otkryto v mnoggobrazii special'nyh
ciklov v sheme Gil'berta (can i prove that it is dense?), v to zhe
vremya lokal'naya svboda - otkrytoe uslovie v $M$, tak chto vse
soglasovano].
\end{for-me}

{\bf Proof.}

    (a) follows from the definition of $\delta(\xi,L+K)$, and duality.
To prove (b), we consider the diagram 
$$
   \smallXyMatrix
   {
      0 \ar[r] & { \Ext^1(\jk(L),\OO) }\ar[r]^{\rm res}\xySimeq[d] & { H^0(X,\sheafExt^1(\jk(L),\OO) ) } \ar[r]\xySimeq[d] & { H^2(X,\sheafHom(\jk(L),\OO)) }\xySimeq[d] \\ 
      0 \ar[r] & { H^1(X,\jk(KL))\dual }\ar[r] & { H^0(\xi,(KL)|_{\xi})\dual }\ar[r] & { H^0(X,KL)\dual} 
   }
$$
from ~\ref{spectral-sequence-for-Ext^1(F,O)}.
Let $W=H^0(\xi-p,(KL)|_{\xi-p})\dual$ and $W'=H^0(p,KL|_{p})\dual$.
We have $H^0(\xi,KL|_{\xi})\dual = W \oplus W'$.  The condition
~\ref{local-freeness-over-a-local-ring} implies that the extension
class $e \in \Ext^1(\jk(L),\OO)$ is not locally free at $p$ iff $e \in
W \subset W \oplus W'$.  Let $K = H^1(X,\jk(KL))\dual$. Since $W$ is of
codimension one in $W \oplus W'$, then either $K \subset W$ and all the
extension classes parametrized by $\Ext^1(\jk(L),\OO)$
are not locally free, or that the set of non-locally free
extensions form a hyperplane in $\Ext^1(\jk(L),\OO)$.
Consider the diagram
$$
  \smallXyMatrix
  {
     0 \ar[r] & {K \cap W}\ar[r]\xyCapA & W\ar[r]\xyCapA & { H^0(X,KL)\dual } \xyRavno[d] \\
     0 \ar[r] & K \ar[r] & W \oplus W' \ar[r] & H^0(X,KL)\dual 
  }
$$
and the dual diagram
$$
  \smallXyMatrix
  {
    {} & 0 & 0 & {} \\
    0  & {(K \cap W)\dual}\ar[u]\ar[l] & {H^0(\xi-p,KL)}\ar[u]\ar[l] & H^0(X,KL)\ar[l] \\
    0  & {K\dual}\ar[l]\ar[u] & {H^0(\xi-p,KL) \oplus H^0(p,KL)}\ar[l]\ar[u]^{{\rm pr}_1} & { H^0(X,KL) }\ar[l]_-{ \mathop{{\rm ev}}_{\xi} }\xyRavno[u] 
  }
$$
where $\mathop{\rm ev}_{\xi}: H^0(X,KL) \to H^0(\xi,KL|_{\xi})$ is
the canonical restriction map ("evaluation").  Then $K \subset W$ iff
$H^0(p,KL) \subset \Im {\rm ev}_{\xi}$, i.e., if there is a curve in
the linear system $|K+L|$ passing through $\xi - p$ but not through
$p$.

\begin{for-me}
\subsection{Regular sequences of length 2.}
\label{regular-sequences-of-length-2}
If $A$ is a commutative ring and $\{f,g\}$ is a regular sequence in $A$,
then the Koszul complex $K_{\cdot}(f,g)$
$$
   0 \to A \strelka{(-g,f)} A^2 \strelka{(f,g)} A \to 0 
$$ 
(located in degrees 2,1,0) is acyclic at $p = 2,1$. Indeed,
$H_2(K_{\cdot}(f,g)) = \mathop{\rm Ann}(f,g) = 0$ and if $(a,b) \in
A^2$ is such that $fa+bg = 0$, then $[b][g] = 0$ in $A/f$, and since
$g$ is not a zero divisor in $A/f$, we have $b=fc$ for some $c \in A$,
which implies $a=-gc$ and $H_1(K(f,g))=0$.

{\bf Lemma.} If $f$ is not a zero divisor in $A$ and $H_1(K_{\cdot}(f,g))=0$,
then $\{f,g\}$ is regular.

{\bf Proof.} If $g \mod fA$ is a zero divisor, then one can
find such a $b \in A$ that $[g][b] = 0$ in $A/f$ while $[b] \ne 0$.
It follows that $gb = -fa$ for some element $a \in A$. Since
$H_1(f,g)=0$, there is a $c \in A$ such that $a=-gc$ and $b=fc$,
which contradicts to $[b] \ne 0$.

Another proof can be obtained from the long exact sequence 16.4 in
~\cite{matsumura-Commutative-ring-theory}.
\end{for-me}

\subsection{$\End_A(J) = 0$}
\label{End(J)=A}
{\bf Lemma.} If $A$ is a ring and $J$ is an ideal such that
$\depth_J(A) \ge 2$, then $\End_A(J) = A$.

Indeed, the exact sequence 
$$ \triple{J}{A}{A/J} $$
implies 
$$ 0 \to \Hom(A,A) \to \Hom(J,A) \to \Ext^1(A/J,A) = 0$$
which implies that a natural map $A \to \Hom(J,A)$ is an isomorphism.
Now the inclusion of $J$ into $A$ induced a monomorphism
$$ 0 \to \Hom(J,J) \to \Hom(J,A) \isom A $$, which is clearly epimorphic
since every element $a \in A$ induced an endomorphism of $J$.

\subsection{ Torsion free extensions $\triple{\OO}{E}{i_* B}$}
\subsubsection{ Affine case}
\label{0-1-is-torsion-free:affine-case}
{\bf 1.} Let $A$ be a commutative ring and $f$ be an element of
$A$ which is not a zero divisor. Then the $A$- module $A/f$
has a free resolution
$$
  \xyTriple{A}{A}{A/f}{f}{}
$$
which gives an isomorphism
$\Ext^1(A/f,A) \simeq A/f$. 

   Let $e \in \Ext^1(A/f,A)$ be a class
of the extension 
\begin{equation}
    \triple{A}{E}{A/f}    
    \label{diagram:ext-E}
\end{equation} 
Let $e = [g]$, as in 
~\ref{extensions-of-modules-over-a-commutative-ring},
where $g \in \Hom_A(A,A) \simeq A$.
Then ~\ref{extensions-of-modules-over-a-commutative-ring}
gives that $E$ can be realized as a cokernel of the map
\begin{equation}
     0 \to A \strelka{(f,g)} A^2 \to E \to 0            
     \label{diagram:res-E}
\end{equation}

Let $J$ be the ideal generated by $(f,g)$ in $A$. The cohomological 
Koszul complex $K^{\cdot}(f,g)$ gives a complex $R(f,g)$
$$
   0 \to A \strelka{(f,g)} A^2 \strelka{(-g,f)} J \to 0
$$
The identity homomorphism $\id: A^2 \to A^2$ descends to the
epimorphic map $w: E \to J$, as in the diagram
$$
  \smallXyMatrix
  {
    {} & {} & {} & 0\ar[d] & {} \\
    {} & {} & {} & {\ker w}\ar[d] & {} \\
    0\ar[r] & A\ar[r]^{(f,g)}\xyRavno[d] & A^2\ar[r]^{\rm can}\xyRavno[d] & E\ar[d]^{w}\ar[r] & 0 \\
    0\ar[r] & A\ar[r]^{(f,g)} & A^2\ar[r]^{(-g,f)} & J\ar[d]\ar[r] & 0 \\
    {} & {} & {} & 0 & {} 
  }
$$

(All the columns and rows of the diagram above are exact except,
perhaps, the bottom row.)

The following lemma can be proved directly:

{\bf Lemma}: 

(1) There is a canonical isomorphism $H^1(K) \isom \ker w$.

(2) If $A$ is integer, then there is a canonical isomorphism 
$H^1(R) \isom \Tors E$.

Proof of (2): If $(a,b) \in \ker(f,g)$, then $fa+gb = 0$ and
$f \cdot [(a,b)] = [(fa,fb)] = [(-gb,fb)] = b \cdot [(-g,f)] = 0 $
in $E$, which gives a map $\ker (f,g) \to \Tors E$ and
an injective map $H^1(R) \to \Tors E$. To prove that it is epimorphic, assume
that $[(a,b)] \in \Tors E$. Then $\alpha \cdot [(a,b)] = 0$ for some
nonzero $\alpha \in A$, i.e., $\alpha a = -g \beta$ and $\alpha b = f \beta$
for some $\beta \in A$, which implies $\alpha (fa+gb) = 0$. If $A$ has no
divisors of zero, this implies $fa+gb = 0$.

{\bf 2.}
Assume that $H^1(R)=0$. Then the isomorphism $w: E \to J$ induces the 
isomorphism of the diagram 
$$
  \xymatrix
  {
    {}&{}&{A}\xyRavno[r]\ar[d]^{(f,g)}&{A}\ar[d]^{f}&{}\\ 
    {0} \ar[r] & 
    {A} \ar[r]^{i_2}\xyRavno[d] &
    {A^2} \ar[r]^{\pi_1}\ar[d]  & 
    {A}\ar[r]\ar[d] & {0} \\ 
    {0}\ar[r] & {A}\ar[r]^{i} & {E} \ar[r]^{p}\ar[d] &
    {A/f} \ar[r]\ar[d]&{0}\\ 
    {}&{}&{0}&{0}&{}
  }
$$
of ~\ref{diagram:extension-of-modules}
with diagram
$$
  \xymatrix
  {
    {}&{}&{A}\xyRavno[r]\ar[d]^{(f,g)}&{A}\ar[d]^{f}&{}\\ 
    {0} \ar[r] & 
    {A} \ar[r]^{i_2}\xyRavno[d] &
    {A^2} \ar[r]^{\pi_1}\ar[d]^{(-g,f)} & 
    {A}\ar[r]\ar[d] & {0} \\ 
    {0}\ar[r] & {A}\ar[r]^{f} & {J} \ar[r]^{\phi}\ar[d] &
    {A/f} \ar[r]\ar[d]&{0}\\ 
    {}&{}&{0}&{0}&{}
  }
$$
where $\phi: J \to A/f$ takes $-ga_1 + fa_2$ to $a_1 \mod fA$.
The map $\phi$ is correctly defined in the case $H^1(R)=0$, since 
$-ga_1 + fa_2 = -ga_1' + fa_2'$ implies $f(a_2-a_2') = g(a_1-a_1')$
and, since $H^1(R) = 0$ is equivalent to $H_1(K(f,g)) = 0$, we have
$(a_1-a_1', a_2-a_2') = (f b, gb)$ for some $b \in A$, which implies
$a_1 \mod fA = a_1' \mod fA$.

{\bf 3.} {\bf Corollary.} If $(f,g)$ is a regular sequence in $A$,
then $w$ is an isomorphism, and the extension class $e$ in
$\Ext^1(A/f,A)$ is equivalent to one of the form
$$ 0 \to A \strelka{f} J \strelka{\phi} A/f \to 0$$.

{\bf Remark.} If $(f,g)$ is regular, then there is an isomorphism
$$
  \xymatrix
  {
     0 \ar[r] & A \ar[r]^{f}\xyRavno[d] & J \ar[r]^{\phi}\xyRavno[d] & A/f \ar[r]\xyIsom[d]^{g} & 0 \\
     0 \ar[r] & A \ar[r]^{f} & J \ar[r]^{\psi} & g(A/f) \ar[r] & 0 
  }
$$
where $\psi$ takes $a \in J$ to $a \mod fA$. (Note that
$H_1(K(f,g)) = 0$ implies isomorphisms 
$$ J/fA \isom Ag/ Ag \intersect Af \isom Ag / Afg \isom g \cdot (A/f) $$

\begin{comment}
If $(f,g)$ is a regular sequence. Let $C$ be the subscheme in $\Spec A$
defined by the ideal $(f)$, and $\xi$ be the subscheme in $\Spec A$
defined by the ideal $(f,g)$. Then the extension
~\eqref{diagram:ext-E} is isomorphic to the exact sequence
$$
  \triple{J_C}{J_{\xi}}{J_{\xi \subset C}}
$$
\end{comment}

{\bf Corollary.} If $A/f$ does not have zero divisors, then all
the nonzero extension classes in $\Ext^1(A/f,A)$ are torsion-free. 

\subsubsection{ Global case.} 
\label{0-1-is-torsion-free:global-case}
If $Y$ is a scheme, $L$ is an invertible sheaf on $Y$ and
$s \in H^0(Y,L)$, we would say that $s$ is $\OO_Y$-regular if
for every affine $U \subset Y$ and every nonzero $f \in H^0(U,\OO)$
$fs \ne 0$. (In particular, it is always the case if $Y$ is integer
and $s \ne 0$.
  
  Let $X$ be a scheme, $C$ be an effective Cartier divisor on $X$
and $B$ be an invertible sheaf on $C$.   Let $i$ be an embedding of $\Supp C$ 
into $X$. Assume that we are given an extension 
\begin{equation}
   \triple{\OO}{E}{i_*(B)}
   \tag{1}
   \label{diagram:O-E-B}
\end{equation}
Let $e$ be its class in $\Ext^1_X(B,\OO)$.

By ~\ref{spectral-sequence-for-Ext^1(F,O)} and 
~\ref{sheafExt^1(B,O)=A}, there are isomorphisms
\begin{equation}
   \Ext^1_X(B,\OO) \simeq 
   H^0(X,\sheafExt^1_X(B,\OO)) \simeq
   H^0(X,i_*(NB^{-1})) \simeq 
   H^0(C, NB^{-1}) 
\end{equation}
where $N = N_{C/X}$ is a normal sheaf to $C$ in $X$. 
Let $s \in H^0(C,NB^{-1})$
be the element corresponding to $e$ under this isomorphism.
Let $L = \OO_X(C)$.

   {\bf Lemma.} If $s$ is $\OO_C$ - regular, then there is an equivalence of
extensions
\begin{equation}
  \xymatrix
  {
     0 \ar[r] & { \OO_X }\ar[r]\xyRavno[d] & 
     E \ar[r]\xyIsom[d] &
     i_*(B) \ar[r]\xyRavno[d] & 0 \\
     0 \ar[r] & { \OO_X }\ar[r]^-{f} & 
     {\jk \tensor L}\ar[r]^-{\phi} & i_*(B) \ar[r] & 0 \\
  }
  \label{diagram:0-1-global-case}
\end{equation}
where $f$ is a canonical section of $L$ vanishing at $C$,
$\phi$ is a some map from $\jk(L)$ to $i_*(B)$, and $\xi$ is
a the set of zeroes of $s$ on $C$.

(1) 

{\bf Proof.}
\providecommand{\Ja}{J_{\alpha}}
\providecommand{\Jb}{J_{\beta}}
\providecommand{\Jab}{J_{\alpha\beta}}
\providecommand{\Xa}{X_{\al}}
\providecommand{\Aa}{A_{\alpha}}
\providecommand{\Ba}{B_{\al}}
\providecommand{\Ea}{E_{\alpha}}
\providecommand{\Ia}{I_{\al}}
\providecommand{\fa}{f_{\al}}
\providecommand{\fb}{f_{\be}}
\providecommand{\ea}{e_{\al}}

   Let us choose an affine cover $U_{\al}$ of $X$, as in the proof of
~\ref{sheafExt^1(B,O)=A}. The restriction of the extension
~\ref{diagram:O-E-B} gives an extension $\triple{\Aa}{\Ea}{\Ba}$ with
some class $e_{\al} \in \Ext^1_{\Aa}(\Ba,\Aa)$ for each $\al$. Each
$\Aa$-module $\Ba$ has a resolution ~\ref{diagram: resolution-of-B},
and we can lift the extension class $e_{\al}$ to a homomorphism $\ga:
M_{1,\al} \to \Aa$, as in
~\ref{extensions-of-modules-over-a-commutative-ring}. Let $s_{\al} \in
\Gamma(C_{\al},NB^{-1})$ be a restriction of $s$ to $\Ua$.  Since
$s_{\al}$ is the image of $e_{\al}$ under the isomorphism
$\sheafExt^1_X(i_*(B),\OO_X) \isom i_*(NB^{-1})$ and since $s$ is $\OO_X$-regular,
$(\fa,\ga)$ is a regular sequence in $\Aa$. Let $\Ja$ be the ideal
generated by $(\fa,\ga)$ in $\Aa$.

   Now ~\ref{0-1-is-torsion-free:affine-case} gives an equivalence of extensions
$$
  \xymatrix
  {
   0 \ar[r] & {\Aa} \ar[r]\xyRavno[d] & {\Ea} \ar[r]\xyIsomBelow[d]^{w_{\al}} & {\Ba} \ar[r]\xyRavno[d] & 0 \\
   0 \ar[r] & {\Aa} \ar[r]^{\fa} & {\Ja} \ar[r]^{\phi_{\al}} & \Ba \ar[r] & 0
  }
$$

  Let us extend all the notations to multiindexes, as in 
~\ref{sheafExt^1(B,O)=A}. 
Since $e_{\be} = \phi_{\ab} \psi_{\ab} e_{\al}$ in $\Aab$, we have
$g_{\be} = \phi_{\ab} \psi_{\ab} g_{\al} + \theta_{\ab} \fa$
for some $\theta_{\ab} \in \Aab$. Since $\psi_{\ab}$ is invertible
in $\Aab$ and $\phi_{\ab}$ is invertible in $\Aab/(\fb)$,  we have
$\Ja \cdot \Aab = \Jb \cdot \Aab$; we denote this ideal as $\Jab$.
  Now the restrictions
of $w_{\al}$ and $w_{\be}$ to $\Uab$ give an isomorphism
\begin{equation}
  \xymatrix
  {
     0 \ar[r] & { \Aab }\ar[r]^-{\fb}\xyRavno[d] & 
     {\Jab}\ar[r]^-{\\phi_{\be}} \xyIsomBelow[d]^{w_{\al}w_{\be}^{-1}} &
     B_{\ab}\ar[r]\xyRavno[d] & 0 \\
     0 \ar[r] & {\Aab}\ar[r]^-{\fb} & 
     {\Jab}\ar[r]^-{\phi_{\al}} & B_{\ab}\ar[r] & 0 \\
  }
  \label{diagram:isomorphism-of-extensions-of-B}
\end{equation}

But by ~\ref{End(J)=A}, the  isomorphism $w_{\al}w_{\be}^{-1}$ 
is a multiplication with some constant $\lambda_{\ab}$, and 
since $\fb = \psi_{\ab} \fa$ we have $\lambda_{\ab} = \psi_{\ab}^{-1}$. 

Let $\xi$ be a subscheme in $X$ given by $(\fa,\ga)$ in $\Ua$.  Since
$\psi_{\ab}^{-1}$ is a gluing cocycle for $L=\OO_X(C)$, we get a
desired equivalence of extensions ~\ref{diagram:0-1-global-case}.

It is clear that $\xi$ as a subscheme of $C$ is given as a set of
zeroes of $s$.

The map $f$ can be described in the following way: $\jk(L) := \jk \tensor L$ can
be realized as a subsheaf of $L$ consisting of sections vanishing at
$\xi$, and since the canonical section of $\OO_X(C)$ vanishes at
$\xi$, it gives a section of $\jk(L)$.

{\bf Remark.}
   Consider $\xi$ as a Cartier divisor on $C$.
There is an exact sequence 
$$
    \triple{ J_{C,X} }{ J_{\xi,X} }{ i_* J_{\xi,C} }
$$ 
which is isomorphic to
$$
    \triple{ \OO_X(-C) }{ J_{\xi,X} }{ i_* \OO_C(-\xi) }
$$ 
Twisting it with $L=\OO_X(C)$, we get an exact sequence
\begin{equation}
    \xyTriple{\OO_X}{J_{\xi}(L)}{i_* ( N \tensor \OO_C(-\xi) ) }{f}{}
    \label{diagram:O-jk(C)-O_C(n-xi)}
\end{equation}
where $f$ is the canonical section of $\OO_X(C)$ vanishing at $C$.
It follows that there is an isomorphism 
$$
   \xymatrix
   {
      0 \ar[r] & \OO_X \ar[r]^{f} \xyRavno[d] & \jk(L) \ar[r]\xyRavno[d] & i_* B \ar[r]\xyIsom[d] & 0 \\
      0 \ar[r] & \OO_X \ar[r]^{f} & \jk(L) \ar[r] & i_* (N \tensor \OO_C(-\xi)) \ar[r] & 0 \\
   }
$$

\subsection{ Locally free extensions $\triple{\triv{W}}{E}{i_* B}$}
\subsubsection{ Affine case}
\label{0-k-is-locally-free:affine-case}
  Let $R$ be a commutative ring and $f \in R$ be a non-zero
divisor. The $R$ - module $R/f$ has a resolution 
$$
  \xyTriple{R}{R}{R/f}{f}{}
$$
and therefore $\Ext^1(R/f,W \tensor\limits_k R)$ is a cohomology of the complex
$$
   \Hom(R, W \tensor R) \strelka{f} \Hom(R, W \tensor R) \to 0
$$
In particular, $\Ext^1(R/f,W \tensor R) \simeq W \tensor_k (R/f)$.

   Let $e \in \Ext^1(R/f,R)$ be the class of an extension 
\begin{equation}
    \triple{R}{E}{R/f}    
    \label{diagram:ext-E-2}
\end{equation} 
and let us choose some lift $g \in \Hom(R, W \tensor R)$ of $e$.
Then ~\ref{extensions-of-modules-over-a-commutative-ring}
gives that $E$ can be realized as a cokernel of the map
\begin{equation}
     0 \to R \strelka{(f,g)} R \oplus (W \tensor\limits_k R) \to E \to 0            
     \label{diagram:res-E-2}
\end{equation}

The exact sequence 
$$
   0 \to \Tor_1(E, k(p)) \to k(p) \strelka{(f(p),g(p))} k(p) \oplus (W \tensor k(p))
$$
associated with a point $p \in \Spec R$ demonstrates that $E$ is
locally free iff $f$ and $g$ do not have common zeroes in $\Spec R$.
\subsubsection{ Global case.} 
\label{0-k-is-locally-free:global-case}
   Let $X$ be a scheme, $C$ be an effective Cartier divisor on $X$ (we
also denote by $C$ the support of $C$), let $i$ be the canonical
embedding of $C$ into $X$, and let $B$ be an invertible sheaf on
$C$. Let $e$ be the class of an extension
\begin{equation}
   \triple{\triv{W}}{E}{i_*(B)}
   \label{diagram:O-E-B-2}
\end{equation}
where $\dim W \ge 2$. 

By 
~\ref{spectral-sequence-for-Ext^1(F,O)} and
~\ref{sheafExt^1(B,O)=A} there is an isomorphism
\begin{multline*}
   \Ext^1_X(i_*(B),W \tensor \OO_X) \simeq 
   W \tensor \Ext^1_X(i_*(B),\OO_X) \simeq \\
   \simeq
   W \tensor H^0(X,\sheafExt^1_X(i_*(B),\OO_X)) \simeq
   W \tensor H^0(C, A) \simeq
   H^0(C, W \tensor_k A) 
\end{multline*}
where $A = NB^{-1}$ and $N = N_{C/X} = \ox(C)|_C$ is the normal sheaf
to $C$ in $X$.

If $s \in H^0(C,W \tensor A)$ is the image of $e$ under this
isomorphism, then ~\ref{0-k-is-locally-free:affine-case} gives that
$E$ is locally free iff $s$ does not vanish on $C$.

   Let $\phi_e$ be the image of $s$ under the isomorphism 
$$ 
    H^0(C, W \tensor A) \isom \Hom(W\dual, H^0(C, A)) 
$$ 
Let $l = \rk \phi $. One can prove that $s$ can be written as $ s =
\sum_1^l w_i \tensor s_i$ for linearly independent $w_1,\dots,w_l \in
W$ and linearly independent $s_1,\dots,s_l \in H^0(C,A)$. We have
$Z(s) = Z(s_1) \intersect \dots \intersect Z(s_l)$.

Note that even if $l \ge 2$ and $s_1, \dots, s_l$ are linearly
independent it is possible that $Z(s_1) \intersect \dots
\intersect Z(s_l)$ is not empty. But if $A$ is globally
generated, then the lemma ~\eqref{generic-section-of-a-gg-lf-sheaf}
implies that this does not happen:

\begin{proposition} If $\dim X = 2$, $C$ is smooth and $A = NB^{-1}$
is globally generated, then the generic section of $W \tensor_k A$
does not vanish on $C$, and the corresponding extension is locally
free.  \end{proposition}

We will also be particularly interested in the case $\dim W = \dim
H^0(C,A)$.  We will show late that if $\phi_e$ is not injective, the
extension $e$ can not be locally free. Assume now that $\phi_E$ is
injective, i.e., is an isomorphism. This defines a unique extension
class $e_{can}$. This extension is locally free iff
$\intersectIn{s \in H^0(C,A)} Z(s) = \emptyset$, i.e., iff $A$
is globally generated.

\subsection{ Criteria for a generic extension to be acceptable.}
\label{criteria-for-a-generic-extension-to-be-acceptable}
{\bf Lemma.} Let $F$ be an acceptable sheaf on a surface $X$. 
We will be considering classes $e \in \Ext^1(F,\triv{W})$ 
parametrizing extensions
\begin{equation}
   \triple{\triv{W}}{E}{F}
   \label{diagram:W-E-F}
\end{equation}
Statement: 

{\bf (i)} If $F$ is a locally free, then all such extensions
are locally free;

{\bf (ii)} If $F = \jk(L)$, where $L$ is an ample line bundle and
$\xi$ is a simple 0-cycle on $X$, and pair
$(\xi,L+K)$ is {\it Caley-Bacharash}, then there exist a finite collection of
hyperplanes $Z \subset \Ext^1(\jk(L),W \tensor \ox)$ such that for every
extension class not in $Z$ the sheaf $E$ is locally free;

\label{1-r-is-acceptable}

{\bf (iii)} If $F=(i_C)_* (B)$, where $i: \Supp C \to X$ is an
(effective) Cartier on $X$ and $B$ is an invertible sheaf on $\Supp C$, then

   (a) If $\dim W = 1$ and $\Supp C$ is reduced and irreducible, then
$E$ is acceptable;

   (b) If $\dim W \ge 2$, $\Supp C$ is smooth and $A = N_{C/X}B^{-1}$
is globally generated, then $E$ given by generic $e \in \Ext^1_X(i_*B,
W \tensor \ox)$ is acceptable.

\label{0-r-is-acceptable}

{\bf Proof.}

{\bf (i)} The exact sequence ~\eqref{diagram:W-E-F} implies
$\Tor_1(E,k(p)) = 0$ for every $p \in X$, which implies {\bf (i)}.

{\bf (ii)} Let $k= \dim W$. In the case $k=1$ this statement was
proved in ~\ref{local-freeness}. If $k > 1$,  we proceed by induction.
% as in ~\ref{stability-of-extensions-by-W}. 

   Let us choose a $k-1$-dimensional subspace $W_{k-1} \subset W$ and
let $W^1=W/W_{k-1}$. Assume that we are given an extension
$$    \triple{\triv{W}}{E}{\jk(L)} $$
with class $e \in \Ext^1(\jk(L),\triv{W})$
The projection $W \to W^1$ induces the diagram
$$
  \displayE{ \triv{W_{k-1}} }{ \triv{W} }{E}{ \jk(L) }{ \triv{W^1} }{ E' }
$$
where $E':= E / (\triv{W_{k-1}})$
The exact sequence
$$
  \xyTriple{W_{k-1}}{W}{W^1}{}{p}
$$
induces the exact sequence
\begin{equation}
  \xyTriple{ %
           \Ext^1(\jk(L),\triv{W_{k-1}}) }{ 
           \Ext^1(\jk(L),\triv{W})       }{
           \Ext^1(\jk(L),\triv{W^1})     }{ }{ 
           p^1
           }
  \label{eq:les-Ext}
\end{equation}           
which is exact on the left and on the right since the
exact sequence
$$
   \triple{\triv{W_{k-1}}}{\triv{W}}{\triv{W^1}}
$$
has a left splitting.  Note that $p^1(e)=e'$, where $e' \in
\Ext^1(\jk(L),\triv{W^1})$ is the class of the bottom row of the
diagram above.  The ~\ref{local-freeness} implies that $E'$ is locally
free when $e'$ is outside of a finite collection of hyperplanes in 
$\Ext^1(\jk(L),\triv{W^1})$. Since local freeness of $E'$ implies
a local freeness of $E$, the exact sequence
~\eqref{eq:les-Ext} 
guarantees that $E$ is locally free for all $e$ outside some
finite collection of hyperplanes in the vector space $\Ext^1(\jk(L),\triv{W})$.

{\bf (iii)} It was proved in ~\ref{0-1-is-torsion-free:global-case}
in the case $\dim W = 1$ and in ~\ref{0-k-is-locally-free:global-case}
in the case $\dim W > 1$.

%------------------------------------------------------------------
\begin{for-me}
\subsection{Remark on the extensions $0 -> k$ and $0 -> 1 -> k$ (for me)}
\subsubsection{ Lemma about induction of extensions.}
% This lemma is used in 
\label{lemma-about-induction}
{\bf Lemma.} Let $$\xyTriple{W_1}{W_2}{W_3}{}{h}$$ be an exact sequence
of sheaves on a scheme $X$, 
assume that $\Ext^1(W_3,W_1)=0$, and let
$B$ be a coherent sheaf on $X$. 
%  such that the induced exact sequence
%$$
%  \xyTriple{\Ext^1(B,W_1)}{\Ext^1(B,W_2)}{\Ext^1(B,W_3)}{}{\ind}
%$$
% is exact
% LISHNEE USLOVIE!
Consider the induced map
$$
    \Ext^1(B,W_2) \strelka{h^1} \Ext^1(B,W_3)
$$
Let $e' \in \Ext^1(B,W_3)$ be a class of an extension 
$$
   \triple{W_3}{G}{B}
$$
and let $E_{e'} = (h^1)^{-1}(e') = \{ e \in \Ext^1(B,W_2): h^1(e) = e' \}$. 

Under these conditions here is a natural isomorphism
$E_{e'} \isom \Ext^1(G,W_1)$.

{\bf Proof.}

  {\bf 1.} Construction of $\alpha: E_{e'} \to \Ext^1(G,W_1)$.
Let $e \in E_{e'}$ be a class of an extension
$$
   \xyTriple{W_2}{F}{B}{i}{\pi}
$$
The map $h: W_2 \to W_3$ induces a diagram
$$
  \smallXyMatrix
  {
    {} & 0\ar[d] & 0\ar[d] & {} & {} \\
    {} & {W_1}\xyRavno[r]\ar[d] & {W_1}\ar[d] & {} & {} \\
    0\ar[r]&{W_2}\ar[r]^{i}\ar[d]& {F}\ar[r]^{\pi}\ar[d]& {B}\ar[r]\xyRavno[d] & 0\\
    0\ar[r] & {W_3}\ar[r]\ar[d] & {F'}\ar[r]\ar[d] & {B}\ar[r] & 0\\
    {} & 0 & 0 & {} & {} 
  } 
$$
Where $F':= F / i(W_1)$. The bottom row of this diagram
is the induced extension. Since $h^1(e') = e$, there is an equivalence
of extensions
$$
   \xyEquivalenceA{W_3}{G}{B}{F'}{i}
$$
The isomorphism $i$ induces an extension
$$
   \xyInducedOnTheRight{W_1}{\tilde F}{G}{F}{F'}{}{i}
$$
where $\tilde F = F \mathop{\oplus}\limits_{F'} G = \ker (F \oplus G \to F')$
This gives an element in $\alpha(e) \in \Ext^1(G,W_1)$.

{\bf 2.} Construction of $\beta: \Ext^1(G,W_1) \to E_{e'}$.
Let $e' \in \Ext^1(B,W_3)$ be a class of an extension
$$
   \xyTriple{W_3}{G}{B}{}{\pi_1}
$$
and assume that we are given an extension class
$$
   \xyTriple{W_1}{F}{G}{}{\pi_2}
$$
in $\Ext^1(G,W_1)$. This data induces a diagram
$$
  \smallXyMatrix
  {
    {} & 0\ar[d] & 0\ar[d] & {} & {} \\
    {} & {W_1}\xyRavno[r]\ar[d] & {W_1}\ar[d] & {} & {} \\
    0\ar[r] & K \ar[r]\ar[d] & F \ar[r]^{\pi_1\pi_2}\ar[d]_{\pi_2} & 
    B \ar[r]\xyRavno[d] & 0\\
    0\ar[r] & {W_3}\ar[r]\ar[d] & G \ar[r]^{\pi_1}\ar[d] & B \ar[r] & 0\\
    {} & 0 & 0 & {} & {} 
  } 
$$
where $K := \ker (\pi_1 \pi_2)$. 
Since $\Ext^1(W_3,W_1) = 0$, there is an isomorphism
$$
   \xyEquivalenceA{W_1}{K}{W_3}{W_2}{i}
$$
which induces the extension
$$
   \xyInducedOnTheLeft{K}{F}{B}{W_2}{\tilde F}{i}{} 
$$
This gives a desired element $e \in E_{e'}$.

It is easy to see that $\alpha \beta = \beta \alpha = 1$.

\subsubsection{  Flag $\xi' \subset \xi \subset C$}
\label{flag-xi-xi'-C}
  Let $X$ be a K3 surface, $C$ be a smooth curve on $X$,   
$\xi' \subset \xi \subset C$, where $\xi'$ is an effective divisor on $C$,
$p \in C$, $p \notin \xi'$, and $\xi = \xi' + p$. Let $L=\OO_X(C)$.
Note that $L|C = K_C$.

  The filtration $J_C \subset J_{\xi} \subset J_{\xi'}$ induces a diagram
$$
   \displayD{J_C}{J_{\xi}}{J_{\xi'}}{\OO_p}{
             \OO_C(-\xi)}{\OO_C(-\xi')}{\OO_C(-\xi')_p }
$$
which after twisting with $L$ gives
$$
   \displayD{\OO_X}{J_{\xi}(L)}{J_{\xi'}(L)}{L_p}{K_C \tensor \OO_C(-\xi)}{
             K_C \tensor \OO_C(-\xi')}{(K_C \tensor \OO_C(-\xi'))_p}
$$
Taking the global sections gives a diagram
\begin{equation}
\smallXyMatrix
{
    {} & 0\ar[d] & 0\ar[d] & {} \\
    {} & k \ar[r]^{\id}\ar[d] & k \ar[d] & {} \\
    0\ar[r] & { H^0( X, J_{\xi}(L) ) }\ar[r]\ar[d] & { H^0( X, J_{\xi'}(L) ) } \ar[r]\ar[d] & { L_p }\xyIsom[d] \\
    0\ar[r] & { H^0(C,K_C \tensor \OO_C(-\xi)) }\ar[r]\ar[d] & { H^0(C,K_C \tensor \OO_C(- \xi')) } \ar[r]\ar[d] & { (K_C \tensor \OO_C(-\xi'))_p } \\
    {} & 0 & 0 & {} 
} 
\label{diagram:xi'<xi<C} 
\end{equation}

\subsubsection{Remark: induction for acceptability, $0-1-k$ vs $0-k$}.

  Assume that $C$ is smooth and irreducible. There is an isomorphism
$$\Ext^1_X(i_*(B),\OO_X) \strelka{\alpha} H^0(C,A),$$ where $A =
N_{C/X} B^{-1}$.  If $e \in \Ext^1_X(i_*(B),\OO_X)$ and $s =
\alpha(e)$, and if $\xi = \xi(e)$ is the set of zeroes of $s$ on $C$,
then ~\ref{...}  gives that $e$ is equivalent to an extension of the form
$$
   0 \to \OO_X \to \jk(C) \to B \to 0
$$

Since $h^0(C,A) > 0$, the 0-cycle $\xi$ on $C$ is special
with respect to the $|K_C|$, or, equivalently, $\xi$ as a 0-cycle
on $X$ is special w.r.t. the linear system $|C|$.

If $k = \dim W > 1$ and if we are given a $k-1$-dimensional vector
subspace $W_{k-1} \subset W$, then we let $W_1=W/W_{k-1}$ and
$E' := E / ( \triv{W_{k-1}} ) $. Consider the diagram
$$
  \displayE{ \triv{W_{k-1}} }{ \triv{W} }{E}{ B }{ \triv{W_1} }{ E' }
$$

Consideration of the diagrams of this kind defines a correspondence
$$
  \xymatrix
  {
    {} & S \ar[ld]_{p_1}\ar[rd]^{p_2} & {} \\
    { \Ext^1_X(B,W \tensor \OO_X) } & {} & T \ar[d]^{p_3} \\
    {} & {} & { \PP\Ext^1_X(B,\OO_X) } 
  }
$$
where $S = \Ext^1_X(B,W \tensor \OO_X) \times \Gr(k-1,W)$ enumerate
diagrams like one above, $p_1$ is a projection on the first summand,
$T$ parametrizes the data of the form $($ one-dimensional vector space
$W_1$, extension class $0 \to W_1 \tensor \OO_X \to E' \to B \to 0$,
and an epimorphism $W \to W_1 \to 0$ $)$. The fibers of $p_3$ are
isomorphic to $\PP(W\dual)$.

Now ~\ref{lemma-about-induction} implies that the fibers of $p_2$
are isomorphic to $W_{k-1} \tensor \Ext^1_X(E',\OO)$.

The lemma (...) gives 
$$
   \delta_X(\xi,L) = h^1(C,B) - p_g = h^0(C,A+K_X\cdot C) - p_g)
$$

(!!!!!!!!!!!!)
[Assume that $\delta_X(\xi,L) > 0$. In particular, the condition
$h^0(C,A) > 0$ sometimes can guarantee that; see the appendix on the
relation of $\delta_C(\xi)$ and $\delta_X(\xi,C)$. For example, this
is the case of K3. Assume, furthermore, that $(\xi,L)$ is Caley-Bacharash.]

Consider the diagram
$$
   \xymatrix
   {
     0 \ar[r] & 
     {\Ext^1_X(B,W_{k-1} \tensor \ox)} \ar[r]\xySimeq[d] &
     {\Ext^1_X(B,W       \tensor \ox)} \ar[r]\xySimeq[d] &
     {\Ext^1_X(B,W_1     \tensor \ox)} \ar[r]\xySimeq[d] &
     0 \\
     0 \ar[r] &
     W_{k-1} \tensor H^0(C,A) \ar[r] &
     W       \tensor H^0(C,A) \ar[r] &
     W_1     \tensor H^0(C,A) \ar[r] &
     0
   }
$$
if $C$ is reduced and irreducible then any extension class
in $\Ext^1_X(B,W_1 \tensor \ox)$ is acceptable.
if $A$ is globally-generated, the {/it generic} extension class
in $\Ext^1_X(B,W \tensor \ox)$ is acceptable. For a given class
$e' \in \Ext^1_X(B,W_1 \tensor \ox)$, its preimage is 
$\Ext^1_X(B,W_1 \tensor \ox)$ which contains a finite collection
of hyperplains so that the complement parametrizes locally extensions.
Therefor the picture is self-consistent.
%
% End of remark on 0-> k vs. 0->1-> k
\end{for-me}
%------------------------------------------------------------------------

\begin{comment}

The exact sequence
$$
  \triple{W_{k-1}}{W}{W_1}
$$
induces the long exact sequence
$$
  \triple{\Ext^1(B,\triv{W_{k-1}})}{\Ext^1(B,\triv{W})}{\Ext^1(B,\triv{W^1})}
$$
which is isomorphic to 
$$
  \triple{ 
           \Hom(W_{k-1}\dual, \Ext^1(B,\OO)) }{
           \Hom(W\dual,       \Ext^1(B,\OO)) }{
           \Hom((W^1)\dual,   \Ext^1(B,\OO)) 
         }
$$
%
Let $e \in \Ext^1(B,\triv{W})$
be a class of some extension
$$
  \triple{\triv{W}}{E}{B}
$$
where $E$ is stable. Then 
~\ref{non-simpleness-for-non-injective-maps}
or
~\ref{stability-of-extensions-by-W}
implies that $e$ maps to a nonzero extension class
$e' \in \Ext^1(B,\triv{W^1})$:
\begin{equation}
   \triple{\triv{W_1}}{E'}{B}
   \tag{$e'$}
\end{equation}

\end{comment}

% -------------------------------------------------------------------------

% LocalWords:  hyperplains
  
          % no [], spell-checked
          % \subsection{(.) Locally free extensions $\triple{\OO}{E}{\jk(L)}$}
          % \subsection{Regular sequences of length 2.} - commented out
          % \subsection{$\End_A(J) = 0$}
          % \subsection{(.) Torsion free extensions $\triple{\OO}{E}{i_* B}$}
          % \subs{(.) Locally free extensions $\triple{\triv{W}}{E}{i_* B}$}
          % \subs{(.) Criteria for a generic extension to be acceptable.}
          % \subs{Remark on the extensions $0 -> k$ and $0 -> 1 -> k$.} -
          % commented out for this version
          % ---------
 
% \input{ structure-of-extensions.part3  }  - excluded
          % \subsection{Lemma about non-simpleness.}  - which i never use
          % -------------------------

\subsection{Stability of extensions and factors}
\label{stability-of-extensions-and-factors}
%
% rasshirenie
%
%-------------------------------------------------------------
\begin{for-me}
     {\bf Lemma.}
     \label{stability-of-extensions-by-O-for-rank>=1}
     %
     % zdes' problema s sluchaem (c_1(E),H) = 1, i ya ego vykinul
     %
      Let $X$ be an (integer) surface, $H$ be a polarization on $X$, and let
     $E'$ be an $H$-stable coherent sheaf on $X$, $\rk E' \ge 1$.
     Assume that we are given a nontrivial extension
     \begin{equation}
        0 \to \OO_X \strelka{i} E \strelka{\pi} E' \to 0
        \label{exact-sequence-E-E'}  % 1
     \end{equation}
     and that either one of the following two conditions holds:
     \begin{enumerate}
         \item[(a)]
                $\Pic S \simeq Z$, $h$ is a positive generator of
     $\Pic S$, and $c_1(E') = h$, or
         \item[(b)]
                     $(c_1(E'),h)=1$. 
     \end{enumerate}
     Then $E$ is $H$-stable.
     
     {\bf Proof.} The map $\pi: E \to E'$ induces a map $\Tors E \to \Tors
     E'$ which is inclusion since $X$ is integer. Since $E'$ is stable of
     rank $\ge 1$, it is torsion free, and therefore $E$ is torsion
     free. Let $r = \rk E$ and $k = \rk F$.
     
     Assume that $E$ is not stable; let $F \subset E$ be a destabilizing
     subsheaf.  The intersection $\OO_X \intersect F$ is a coherent
     subsheaf of $\OO_X$ and thus is either 0, or a sheaf of ideals of a
     proper subscheme $Z \subset X$.
     
        Consider first the case $\OO_X \intersect F = 0$. In this case $F'
     = \pi(F)$ is a subsheaf of $E'$ isomorphic to $F$. Since $E'$ is
     stable, we have
     \begin{gather}
        \frac{P(F)}{k} \ge \frac{P(E)}{r}, \\ 
        \frac{P(F')}{k} <  \frac{P(E')}{r-1} 
        \label{equation:condition-on-F}
     \end{gather}
     The reduced Poincare polynomial $P(E)/r$ of a torsion-free coherent sheaf
     of rank $r$ on $X$ can be written as
     $$
        \frac{P(E)}{r} = 
        H^2 \frac{m(m+1)}{2} + \frac{\chi(E_H)}{r} m + \frac{\chi(E)}{r}
     $$
     where $E_H$ is a restriction of $E$ to a smooth curve $H$ in the
     linear system $|h|$ (which we can assumed to be very ample), and
     therefore ~\ref{equation:condition-on-F} is equivalent to
          \begin{gather}
             \frac{\chi(F_H)}{k} m + \frac{\chi(F)}{k} \ge \frac{\chi(E_H)}{r} m + \frac{\chi(E)}{r} , \\
             \frac{\chi(F'_H)}{k} m + \frac{\chi(F)}{k} <   \frac{\chi(E'_H)}{r-1} m + \frac{\chi(E')}{r-1}
        \label{equation:condition-on-F-2}
     \end{gather}
     
     The Riemann-Roch theorem for a torsion-free coherent sheaf $G$ of rank
     $r$ on a smooth curve $C$ of genus $g$ gives
     $$
       \chi(C,G) = r(1-g) + c_1(G)
     $$
     and therefore
     ~\ref{equation:condition-on-F-2}
     is equivalent to 
     \begin{gather}
        \frac{(c_1(F),H)}{k} m + \frac{\chi(F)}{k} \ge \frac{(c_1(E),H)}{r} m + \frac{\chi(E)}{r},\\
        \frac{(c_1(F),H)}{k} m + \frac{\chi(F)}{k} <   \frac{(c_1(E),H)}{r-1} m +
                                                                \frac{\chi(E')}{r-1}
        \label{equation:condition-on-F-3}
     \end{gather}
     
     In particular, we have
     \begin{gather*}
        \frac{(c_1(F),H)}{k} \ge \frac{(c_1(E),H)}{r}  ,\\
        \frac{(c_1(F),H)}{k} \le \frac{(c_1(E),H)}{r-1} 
     \end{gather*}
     
     If $\Pic S \simeq \ZZ h$, $c_1(E) = h$, $c_1(F) = f \cdot h$, as in
     the first statement of the lemma, then the condition above is
     equivalent to the condition
     $$
         \frac{1}{r} \le \frac{f}{k} \le \frac{1}{r-1}
     $$
     Since $k \in [1,\dots,r-1]$, this equality can be satisfied only if $k
     = r-1$ and $f=1$. (Drawing a picture in $\ZZ^2$ assigning to a
     fraction $a/b \in \QQ$ a point with coordinates $(a,b)$ helps to
     visualize the condition.)
     
     If $(c_1(E),H) = 1$, as in the second statement of the lemma, then the
     condition above is equivalent to 
     $$
         \frac{1}{r} \le \frac{(c_1(F),H)}{k} \le \frac{1}{r-1}
     $$
     which can be satisfied only if $k = r-1$ and $(c_1(F),H)=1$.
     
     In both cases we must have $\rk F' = \rk E' = r-1$, and therefore $\rk
     E'/F' = 0$.
     
     Let $Q = E / F$ and $Q' = E'/F'$, and consider the diagram
     \begin{equation}
       \xymatrix
       {
          {} & {} & 0 \ar[d] & 0 \ar[d] & {} \\
          {} & {} & F \xyIsom[r]\ar[d] & {F'}\ar[d] & {} \\
          0 \ar[r] & {\OO_X} \xyIsom[d]\ar[r] & E \ar[r]\ar[d] & {E'} \ar[d]\ar[r] & 0 \\
          0 \ar[r] & {\OO_X} \ar[r] & Q \ar[r]\ar[d] & {Q'} \ar[d]\ar[r] & 0 \\
          {} & {} & 0 & 0 & {}
       }
       \label{diagram:narushenie-stabilnosti-rasshireniya}
     \end{equation}
     
     Note that in both cases we have $(c_1(Q'),H) = (c_1(E),H) -
     (c_1(F),H) = 0$.
     
       Since $\rk Q' = 0$, the $\Supp Q'$ is either 0-dimensional or
     1-dimensional in $X$, But since $H$ is ample and $(c_1(Q'),H) =0$,
     $\Supp Q'$ can not be 1-dimensional. (In the case $\Pic S \simeq \ZZ$
     this is obvious. In the case $(c_1(E),H) = 1$ we should exclude the
     case $\Supp Q'$ is a rational curve contracted by $h$.)

     It follows that $c_1(Q') = 0$, and the equations
     ~\ref{equation:condition-on-F-3} are satisfied automatically.
     
     Consider the exact sequence
     $$
        \Ext^1(Q',\OO_X) \strelka{d_1} \Ext^1(E',\OO_X) \strelka{d_2} \Ext^1(F',\OO_X)
     $$ 
     Let $e \in \Ext^1(E',\OO_X)$ be the class of the extension
     ~\ref{exact-sequence-E-E'}. We have $d_2(e) = 0$ since the extension
     $e$ is split over $F' \subset E'$, and therefore $e=d_1(e')$ for some
     class $e' \in \Ext^1(Q',\OO_X)$. (One can check that $e'$ is the
     extension class of the bottom row of the diagram
     ~\ref{diagram:narushenie-stabilnosti-rasshireniya}). But since $\dim
     \Supp Q' = 0$, we have $\Ext^1(Q',\OO_X) = 0$, and therefore $e=0$,
     which contradict to the assumption that the extension
     ~\ref{exact-sequence-E-E'} is not trivial.
     
      Remark: If we know additionally that $E'$ is locally free, then
     we may assume that $F$ is locally free. In this case $F'$ is a locally
     free subsheaf of $E'$, $\rk F' = \rk E' = r-1$, and locally on $X$ and
     inclusion $F' \to E'$ is isomorphic to an inclusion $\OO_U^{r-1} \to
     \OO_U^{r-1}$ which can not degenerate in codimension two, which
     immediately leads to a contradiction  (if we can prove that it does not
     degenerate along a rational curve)
     
     % *********
     % one can assume that Q is torsion free; (maximal destabil. sheaf);
     % then the bottom line is 
     % 0 -> O -> J_x(D) -> O_D(D-x) -> 0
     %    if Pic = Z, then we get c_1(Q) = 0, i.e,
     % [D] = 0, and supp Q' is 0-dimensional, etc.
     % but if (c1(E),H) = 1, we merely get (D,H) = 0 ,
     % i.e., D is a pre-image of a singularity, Q' is supported on it,
     % and Ext^1(Q',O) ~ H^1(D,Q') = H^1(D, D|_D -x) = H^0(D,x) > 0,
     % and our proof does not work.   
     % *********
     %
\end{for-me}
% ----------------------------------------------------------
     
\begin{thing}
{\bf Lemma.}
\label{stability-of-extensions-by-O-for-rank>=1}
 Let $X$ be an (integer) surface, $H$ be a polarization on $X$, and let
$E'$ be an $H$-stable coherent sheaf on $X$, $\rk E' \ge 1$.
Assume that we are given a nontrivial extension
\begin{equation}
   0 \to \OO_X \strelka{i} E \strelka{\pi} E' \to 0
%  \label{exact-sequence-E-E'} % n 2
\end{equation}
and assume further that $\Pic S \simeq \ZZ$, $h$ is a positive generator of
$\Pic S$.

Then $E$ is $H$-stable.
\end{thing}

{\bf Proof.} The map $\pi: E \to E'$ induces a map $\Tors E \to \Tors
E'$ which is inclusion since $X$ is integer. Since $E'$ is stable of
rank $\ge 1$, it is torsion free, and therefore $E$ is torsion
free. Let $r = \rk E$ and $k = \rk F$.

Assume that $E$ is not stable; let $F \subset E$ be a destabilizing
subsheaf.  The intersection $\OO_X \intersect F$ is a coherent
subsheaf of $\OO_X$ and thus is either 0, or a sheaf of ideals of a
proper subscheme $Z \subset X$.

   Consider first the case $\OO_X \intersect F = 0$. In this case $F'
= \pi(F)$ is a subsheaf of $E'$ isomorphic to $F$. Since $E'$ is
stable, we have
\begin{gather}
   \frac{P(F)}{k} \ge \frac{P(E)}{r}, \\ 
   \frac{P(F')}{k} <  \frac{P(E')}{r-1} 
   \label{equation:condition-on-F}
\end{gather}
\begin{for-me}
% v etoj chasti dokazatel'stva dostatochno slopes,
% poetomu to, chto dal'she, mozhno opustit'

The reduced Poincare polynomial $P(E)/r$ of a torsion-free coherent sheaf
of rank $r$ on $X$ can be written as
$$
   \frac{P(E)}{r} = 
   H^2 \frac{m(m+1)}{2} + \frac{\chi(E_H)}{r} m + \frac{\chi(E)}{r}
$$
where $E_H$ is a restriction of $E$ to a smooth curve $H$ in the
linear system $|h|$ (which we can assumed to be very ample), and
therefore ~\ref{equation:condition-on-F} is equivalent to
\begin{gather}
   \frac{\chi(F_H)}{k} m + \frac{\chi(F)}{k} \ge \frac{\chi(E_H)}{r} m + \frac{\chi(E)}{r} , \\
   \frac{\chi(F'_H)}{k} m + \frac{\chi(F)}{k} <   \frac{\chi(E'_H)}{r-1} m + \frac{\chi(E')}{r-1}
   \label{equation:condition-on-F-2}
\end{gather}

The Riemann-Roch theorem for a torsion-free coherent sheaf $G$ of rank
$r$ on a smooth curve $C$ of genus $g$ gives
$$
  \chi(C,G) = r(1-g) + c_1(G)
$$
and therefore
~\ref{equation:condition-on-F-2}
is equivalent to 
\begin{gather}
   \frac{(c_1(F),H)}{k} m + \frac{\chi(F)}{k} \ge \frac{(c_1(E),H)}{r} m + \frac{\chi(E)}{r},\\
   \frac{(c_1(F),H)}{k} m + \frac{\chi(F)}{k} <   \frac{(c_1(E),H)}{r-1} m +
                                                           \frac{\chi(E')}{r-1}
   \label{equation:condition-on-F-3}
\end{gather}
\end{for-me}

In particular, we have the inequalities of slopes
\begin{gather*}
   \frac{(c_1(F),H)}{k} \ge \frac{(c_1(E),H)}{r}  ,\\
   \frac{(c_1(F),H)}{k} \le \frac{(c_1(E),H)}{r-1} 
\end{gather*}

which is equivalent to 
$$
    \frac{1}{r} \le \frac{f}{k} \le \frac{1}{r-1}
$$
Since $k \in [1,\dots,r-1]$, this equality can be satisfied only if $k
= r-1$ and $f=1$. (One can draw here a picture in $\ZZ^2$ to visualize
the condition.)

Let $Q = E / F$ and $Q' = E'/F'$. We have $\rk Q = 1$, $\rk Q' =
0$. Consider the diagram
\begin{equation}
  \xymatrix
  {
     {} & {} & 0 \ar[d] & 0 \ar[d] & {} \\
     {} & {} & F \xyIsom[r]\ar[d] & {F'}\ar[d] & {} \\
     0 \ar[r] & {\OO_X} \xyIsom[d]\ar[r] & E \ar[r]\ar[d] & {E'} \ar[d]\ar[r] & 0 \\
     0 \ar[r] & {\OO_X} \ar[r] & Q \ar[r]\ar[d] & {Q'} \ar[d]\ar[r] & 0 \\
     {} & {} & 0 & 0 & {}
  }
  \label{diagram:narushenie-stabilnosti-rasshireniya}
\end{equation}

Note that $c_1(Q') = c_1(E) - c_1(F) = 0$. By choosing maximal
destabilizing subsheaf we may assume that $Q$ is torsion free, which
implies that it is of the form $J_{\xi}$ for some 0-subscheme $\xi$ on
$X$ which has no sections, and thus in this case lemma is proved.

{\bf Remark.} If $E'$ is locally free, then we may assume that $F$ is
locally free. In this case $F'$ is a locally free subsheaf of $E'$ and
they have the same rank, while the factor-sheaf is supported in
codimension two, which is not possible, and we get another proof.

% Zamechanie: uslovie $c_1(E) = h$ ili $(c_1(E),H) = 1$ dejstvitel'no 
% neobhodimo.
% V samom dele, pust', naprimer,
% r = 4 and c_1(E)H = 5. Then (f,k) = (3,2) or (4,3) oba udovletvoryayut
% chislennomu kriteriyu naverhu. Drugimi slovami,
% 5/3 > 3/2 > 5/4, ili 5/3 > 4/3 > 5/4. 
% (zdes' horosho narisovat' treugol'nik v ploskosti r,d. u nego est' dve
% celyh tochki. no esli ego vysota = 1, to, konechno, net celyh tochek.)
%
% Ili, v sluchae Pic S = Z, c_1(E) = eh, e>=2, nevozmozhno garantirovat'
% stabil'nost' $E$, dazhe esli potrebovat' $(e,r)=(e,r-1)=1$, t.k.
% vozmozhno, chto c_1(F) = fh, i e/r <= f/k <= e/(r-1). 
% (Naprimer, e=3, r=5, f=2,k=3). 3/5 < 2/3 < 3/4.
%
% Bolee togo, sluchaj $Pic S = Z$, c_1(E) = e h, e>=2 nichem
% ne prosche samogo obschego sluchaya (eto legko uvidet').

\begin{for-me}
{\bf Remark.} In the case $\Pic S = \ZZ h$,  $c_1(E) = e h$,
$e \ge 2$ there is no way to ensure numerically the stability of $E$,
as one can see from the inequalities like  $3/5 < 2/3 < 3/4$,
and this case is not simpler then the case of an arbitrary $\Pic S$.
$
\end{for-me}

% V samom dele, v samom obschem sluchae need to ensure
%  (c_1(E),h) / r <=  (c_1(F),h) / k <= (c_1(E),h) / r-1; dazhe esli
% podelit' na obschij znamenatel', etot sluchaj nichem ne prosche:
% g:= g.c.d. ( (c_1(E),h), (c_1(F),h);
% e:= (c_1(E),h) /g,  f:= (c_1(F),h) /g,  
% e/r <= f/k <= e/r-1 net sposoba obespechit', esli e >= 1.
%

Case 2: consider now the case $\OO_X \intersect F = J$ for a sheaf of
ideals $J \subset \OO_X$ of a subscheme $Z \subset X$.

Consider the diagram
$$
  \xymatrix
  {
     0 \ar[r] & {J_Z}   \ar[r]\ar[d] & F \ar[r]\ar[d] &  F' \ar[r]\ar[d] & 0   \\
     0 \ar[r] & {\OO_X} \ar[r]       & E \ar[r] &  E' \ar[r] & 0 
  }
$$

Let $k = \rk F$. The condition that $F$ destabilizes $E$ and
$F'$ does not destabilize  $E'$,
\begin{gather*}
  \frac{P(F)}{k} \ge \frac{P(E)}{r} \\
  \frac{P(F')}{k-1} < \frac{P(E')}{r-1} 
\end{gather*}
is equivalent to the condition
\begin{gather*}
  P(F) \ge \frac{k}{r} P(E) \\
  P(F) <   P(J_Z) + \frac{k-1}{r-1} (P(E) - P(\OO_X))
\end{gather*}
In particular, we should have
$$
 P(J_Z) + \frac{k-1}{r-1} (P(E) - P(\OO_X)) > \frac{P(E)}{r}
$$
or
$$
   r(r-1)P(J) + (k-r) P(E) -r(k-1) P(\OO_S) > 0
$$

Using the formula
$$
   \frac{P(F)}{k} = H^2 \frac{m(m+1)}{2} + \frac{\chi(F|_H)}{k} m 
                        + \frac{\chi(F)}{k}
$$
we see that the coefficients at $m^2$ of the polynomial 
above is zero, while the condition on the coefficient at $m$ gives
an inequality
$$
  (r-k) (c_1(E),H)_S + r(r-1) \deg_H (Z) \le 0
$$
Since $(c_1(E),H)_S > 0$ and $\deg_H (Z) \ge 0$, we should have
$r=k$ and $\deg_H (Z) = 0$.

It follows that the stability condition ~\ref{equation:condition-on-F}
is equivalent to 
\begin{gather*}
  P(F) \ge P(E) \\
  P(F) <   P(J_Z) + P(E) - P(\OO_S)
\end{gather*}
which implies $P(E) < P(J_Z) + P(E) - P(\OO_S)$, or $P(\OO_S) <
P(J_Z)$, which can not be satisfied, since 
$$
    P(\OO_Z) = \chi(\OO_Z) \ge 0
$$

\begin{thing}
%
% faktorizaciya
\label{stability-of-O-factors}
{\bf Lemma.} Let $X$ be a surface, $H$ be a polarization on $X$,
and $E$ be an $H$-stable coherent sheaf on $X$ such that $\rk E \ge 2$. 
Assume that we are given an exact sequence
\begin{equation}
   0 \to \OO_X \strelka{i} E \strelka{\pi} E' \to 0
%  \label{exact-sequence-E-E'} % 3
\end{equation}
where $E'$ is torsion-free, and that $\Pic S \simeq Zh$, and $c_1(E) = h$.
Then $E'$ is $H$-stable.
\end{thing}

{\bf Proof.}  Assume that $E'$ is not stable; then it admits a
destabilizing quotient sheaf $Q$. Let $r = \rk E$ and $k = \rk Q$.
Since $E'$ is torsion free, $\rk \ker (E' \to Q) \ge 1$, and therefore
$k \le r-2$.

   Since $E$ is stable, $Q$ should not destabilize $E$. This gives
a numerical condition
\begin{gather*}  
   \frac{P(Q)}{k} \le \frac{P(E')}{r-1} \\
   \frac{P(Q)}{k} >   \frac{P(E)}{r}
\end{gather*}  
or
\begin{gather*}  
   \frac{\chi(Q_H)}{k} m + \frac{\chi(Q)}{k} \le \frac{\chi(E'_H)}{r-1} m + \frac{\chi(E'_H)}{r-1} \\
   \frac{\chi(Q_H)}{k} m + \frac{\chi(Q)}{k} >   \frac{\chi(E_H)}{r} m + \frac{\chi(E)}{r} \\
\end{gather*}  
Comparing the coefficients at $m$, we get 
$$
    \frac{(c_1(E),H)}{r} \le \frac{(c_1(Q),H)}{k} \le \frac{(c_1(E),H)}{r-1} 
$$
Now if $c_1(Q) = qh$ this is equivalent to the inequality
$$
    \frac{1}{r} \le \frac{q}{k} \le \frac{1}{r-1} 
$$
which is impossible since $k \le r-2$.

\begin{for-me}
{\bf Remark.} In the case $\Pic S \simeq \ZZ h$, $c_1(E) = e h$, $e
\ge 2$, just as in the general case, there is no numerical way to
ensure the stability of $E'$.
\end{for-me}

% stabil'nost' rasshirenij s W tensor O
\begin{thing}
{\bf Lemma.}
\label{stability-of-extensions-by-W-tensor-O}
   Let $X$ be a surface, $H$ be a polarization on $X$,
and $F$ be an $H$-stable acceptable coherent sheaf on $X$. Assume that we
are given an extension
\begin{equation}
   0 \to W \tensor \OO_X \to E \to F \to 0
   \label{exact-sequence-E-F}
\end{equation}
Let $e \in \Ext^1(F,W \tensor \OO_X)$
be the class of this extension
and let $\alpha_e$ be the image of $e$ under the isomorphism
$\Ext^1(F,\triv{W}) \isom \Hom(W\dual,\Ext^1(F,\OO))$.

Then 

  (1) If $\alpha_e$ is not monomorphic, then $E$ is not $H$-stable.

  (2) If $\rk F \ge 1$, $X$ is irreducible and integer, $\alpha_e$ is
monomorphic, $\Pic X \simeq \ZZ h$, and $c_1(E) = h$, then $E$ is $H$
- stable.

  (3) If $\rk F = 0$, $F = i_*(B)$, where $i: C \to X$ is an embedding
of the support of an integer Cartier divisor $C$ on $X$, and $B$ is an
invertible sheaf on $C$, $X$ is irreducible and integer, $\alpha_e$ is
monomorphic, and $\Pic X \simeq \ZZ [C]$, then $E$ is $H$ - stable.
\end{thing}

{\bf Proof.} Assume first that $\dim W = 1$. If $\rk F \ge 1$, then
the first statement of the lemma is trivial and the second one was
proved in ~\ref{stability-of-extensions-by-O-for-rank>=1} If $\rk F =
0$, then the first statement is trivial and the third one follows from
~\ref{0-1-is-torsion-free:global-case}. 

   Assume now that $\dim W > 1$. Let $k = \dim W$. Let us choose a
$k-1$ -dimensional vector subspace $W_{k-1} \subset W$ and let $W_1 = W/W_{k-1}$.
Then $W_{k-1} \tensor \OO_X$ is a subsheaf of $W \tensor \OO_X$;
let $E' = E / (W_{k-1} \tensor \OO_X)$. The canonical projection
map $\pi: W \to W_1 \to 0$ induces the diagram
\begin{equation}
   \displayE{W_{k-1} \tensor \OO}{W \tensor \OO}{E}{F}{W_1 \tensor \OO}{E'}
   \label{diagram:E,F,E'}    
\end{equation}

The map $\pi$ also induces the diagram
$$
  \xymatrix
  {
    \Ext^1(F, W \tensor \OO) \ar[r]^{d}\xySimeq[d] & \Ext^1(F, W_1 \tensor \OO)\xySimeq[d] \\
    \Hom(W\dual,\Ext^1(F,\OO)) \ar[r]^{d'} & \Hom(W_1\dual,\Ext^1(F,\OO))
  }
$$

It is clear that $e' = d(e)$ is the class of the bottom row
of the diagram ~\eqref{diagram:E,F,E'}. 

If $\alpha_e: W\dual \to \Ext^1(F,\OO)$ is not injective,
then we can choose $\pi: W \to W_1 \to 0$ in such a way that
$W_1\dual \subset \ker(\alpha_e: W\dual \to \Ext^1(F,\OO))$.
In this case $d'(\alpha_e) = 0$, or, in other words, the bottom
row of ~\eqref{diagram:E,F,E'} is a trivial extension
and $E$ admits $\OO_X$ as a factor-sheaf. It follows immediately
that $E$ can not be stable, since it admits $\OO_X$ as its subsheaf
and its factor-sheaf. (Moreover, it is easy to see that under some
mild assumptions on $F$ $E$ can not be semi-stable.) This proves (1).

Assume now that $\alpha_e$ is injective.  The bottom row of the diagram  
induces the exact sequence
$$
  0 \to 
  \Hom(W_1 \tensor \OO,W_1 \tensor \OO) \strelka{\delta} 
  \Ext^1(F,W_1 \tensor \OO) \to 
  \Ext^1(F',W_1 \tensor \OO)
$$
which is injective on the left since 
$\Hom(W_1 \tensor \OO,W_1 \tensor \OO) \simeq H^0(X,\OO) = k$
and since $\delta$ takes the canonical element of 
$\Hom(W_1 \tensor \OO,W_1 \tensor \OO)$ to the nontrivial class $e'$.

Let $e'' \in \Ext^1(E',W_{k-1} \tensor \OO_X)$ be the class of the middle
column of ~\eqref{diagram:E,F,E'} and let $\alpha_{e''}$ be the
corresponding element in $\Hom(W_{k-1}\dual,\Ext^1(E',\OO_X))$.
Consider the diagram
\begin{equation}
  \xymatrix
  {
    0 \ar[r] & W_1\dual \ar[r] & \Ext^1(F,\OO) \ar[r] & \Ext^1(E',\OO) \\
    0 \ar[r] & W_1\dual \ar[r]\xyRavno[u] & W\dual \ar[r]\ar[u]_{\alpha_e} & {W_{k-1}}\dual \ar[r]\ar[u]_{\alpha_{e''}} & 0 \\
    {} & {} & 0 \ar[u] & {} & {} 
  }
\end{equation}

The top row of this diagram is induced by applying $\RHom(\cdot,\OO)$
to the bottom row of the ~\eqref{diagram:E,F,E'} (its injectivity on
the left follows from the discussion above). The fact that this
diagram commutes can be verified directly. This diagram implies that
$\alpha_{e''}$ is injective.

We can now apply induction to the pair $(W,W_{k-1})$. Indeed, in the
case $\rk F \ge 1$ the case $k=1$ considered above implies that $E'$
is stable.  We can now apply induction to the extension
$$
   0 \to W_{k-1} \tensor \OO_X \to E \to E' \to 0
$$
of the middle column of the diagram ~\eqref{diagram:E,F,E'}.

In the case $\rk F = 0$, $F = i_*(B)$ the property that $C$ is integer
and the case $k=1$ considered above imply that $E'$ is stable of rank
one (which is in this case equivalent to the fact that $E'$ is
torsion-free). Now we can apply the induction step to the middle
column extension to deduce that $E$ is stable.

\begin{thing}
\label{stability-of-WO-factors}
{\bf Lemma.} Let $(X,H)$ be a polarized surface such
that $\Pic S \simeq \ZZ [H]$ and $E$ be an $H$-stable coherent sheaf
on $X$ such that $c_1(E) = [H]$. Assume that we are given an exact
sequence 
\begin{equation}
   0 \to W \tensor \OO_X \strelka{i} E \strelka{\pi} F \to 0
\end{equation}
where $F$ is torsion-free. Then $F$ is $H$-stable.
\end{thing}

{\bf Proof.}  Note that the assumption implies $\rk E \ge 2$. 
If $\dim W = 1$, the lemma was proved above. Assume that $\dim W >1$
and let us choose an 1-dimensional vector subspace $W_1 \subset W$.
Let $E' = E /(W_1 \tensor \OO_X)$.
consider the diagram
$$
  \displayE{\triv{W_1}}{\triv{W}}{E}{F}{\triv{(W/W_1)}}{E'}
$$
Since $\OO_X$ is integer, the bottom row of this diagram gives an
embedding $0 \to \Tors E' \to \Tors F$, and since $F$ is torsion free,
$E'$ is also torsion-free. The case $\dim W = 1$ considered above
applied to the middle column of the diagram above implies that 
$E'$ is stable,  and we can now apply induction to the bottom row
of the diagram above.
 
          % no extra [], spell-checked
          %\subsection{(*)Stability of extensions and factors}
          %--------

\subsection{ Flag $\xi \subset C \subset X$} \label{0-1-aj} 
  Let $X$ be a surface, $C$ be an effective Cartier divisor on $X$,
and $\xi$ be an effective Cartier divisor on $\Supp C$. We will also
write $C$ for $\Supp C$ and $\xi$ for $\Supp \xi$. Let $i: C \to X$
and $j: \xi \to C$ be the canonical embeddings.

  The filtration $J_C \subset J_{\xi} \subset \OO_X$ induces the diagram
\begin{equation}
  \xymatrix@R=7pt@C=7pt
    {
      0 \ar[r] & {\jk / J_C} \ar[r]\xySimeq[d] & {\OO_X/J_C} \ar[r]\xyIsom[d] &              {\OO_X/\jk} \ar[r]\xySimeq[d] & 0 \\ 
      0 \ar[r] & i_* J_{\xi,C} \ar[r] & i_* \OO_C \ar[r] & i_* j_* \OO_{\xi} \ar[r] & 0
    }
\end{equation}

and the diagram
$$
  \displayE{J_C}{J_{\xi}}{\OO_X}{i_*j_*\OO_{\xi}}{i_*J_{\xi,C}}{i_*\OO_C}
$$

Let $L = \OO_X(C)$ and $l=i^*(L)$. Twisting the diagram above with $L$,
we get 
\begin{equation}
  \displayD{\ox}{\jk(L)}{L}{(i_*j_* \ok)(L)}{i_* l(-\xi)}{i_*l}{i_*((j_*\ok)(l))}
\end{equation}

\subsection{ Globally generated sheaves and the correspondence}
\label{aj-preserves-gg}
  {\bf Lemma} Let $X$ be a scheme with $h^0(\OO_X)=1$
and 
$$\triple{\OO_X}{E}{E'}$$
be an exact sequence of coherent sheaves. 
If $E$ is globally generated, then $E'$ is globally generated.
If $E'$ is globally generated and $H^1(X,\OO_X)=0$, then $E$ is globally generated.

The proof is straightforward.
\begin{for-me}
{\bf Proof.} Consider the diagram
$$
  \xymatrix@R=8mm@C=8mm
  {
     0\ar[r] & H^0(X,\OO_X) \tensor \OO_X \ar[r]\xyIsom[d] & { H^0(X,E)  \tensor \OO_X } \ar[d]^{e_1}\ar[r] & { H^0(X,E') \tensor \OO_X } \ar[d]^{e_2}\ar[r] & H^1(X,\OO_X) \tensor \OO_X \\
     0\ar[r] & \OO_X \ar[r] &  E \ar[d]\ar[r] & {E'} \ar[r]\ar[d] & 0 \\
     {} & {} & \coker e_1 \ar[r]^u\ar[d] & \coker e_2 \ar[d] & {} \\
     {} & {} & 0 & 0 & {} \\
  } 
$$
It is easy to see that $u$ is epimorphic and in the case
$H^1(X,\OO_X) = 0$ $u$ is an isomorphism. The lemma follows.
\end{for-me}
 
          % no extra [], spell-checked
          %\subsection{(.) Flag $\xi \subset C \subset X$} \label{0-1-aj} 
          %\subsection{(.) Globally generated sheaves and the correspondence}
          % -------

\subsection{Special 0-cycles on surfaces}
\providecommand{\Base}{\mathop{\rm Base}}  % base points of |L|

\subsubsection{Speciality index of a 0-cycle.}

  Let $\xi$ be a 0-dimensional subscheme of a surface $X$ of length
  $d$.

{\bf Definition.} The speciality index of $\xi$ with respect to $L$,
$\delta(\xi,L)$, is defined by
$$
   h^0(X,\jk(L)) =  h^0(X,L) - d + \delta(\xi,L)
$$

It follows that $\delta(\xi,L) \ge 0$ and $\delta(\xi,L) > 0$ iff
$\xi$ does not impose $d$ independent conditions on the sections of
$L$.

{\bf Definition.} $\xi$ is said to be special with respect to $L$ if
$\delta(\xi,L) > 0$.

{\bf Examples.} We assume that the linear system $|L|$ is not empty.

1. If $\deg \xi = 1$, $\xi = [p], p \in X$, then $(\xi,L)$ is special
iff $p$ is a base point for $|L|$.

2. If $\deg \xi = 2$, $\xi = [p] + [q]$, and $|L|$ is base-point free,
then $(\xi,L)$ is special iff $\phi_L(p) = \phi_L(q)$, where $\phi_L:
X \to |L|\dual$ is the morphism to a projective space associated with
$|L|$.

3. Let $\deg \xi = 3$, $\xi = [p]+[q]+[r]$, and let $|L|$ be base-point
free. We denote by $<p_1,\dots,p_n>$ a linear span of a set of
distinct points in projective space. Then $(\xi,L)$ is special iff
either $\phi_L(p) = \phi_L(q) = \phi_L(r)$, or $\dim <p,q,r>_L = 1$,
where $<p,q,r>_L = <\phi_L(p),\phi_L(q),\phi_L(r)>$.

4.  If $\deg \xi = d$ and $\xi$ is a sum of $d$ distinct points on
$X$, and $|L|$ is base-point free, then $\xi$ is special iff $\dim
<\xi>_L < d-1$. 

   Note that in general $\dim <\xi>_L = d-1 - \delta(\xi,L)$.

{\bf Lemma.} $\delta(\xi,L) = h^1(\jk(L)) - h^1(L)$. In particular, if
$X$ is smooth and $L$ = $M \tensor \Omega^2(X)$ for an ample $M$, then
$\delta(\xi,L) = h^1(\jk(L))$.

The proof follows from the long exact sequence associated with
$$
   \triple{\jk(L)}{L}{\ok(L)}
$$
and the Kodaira vanishing theorem.

\subsubsection{Caley-Bacharash 0-cycles}

Let $\xi$ be a simple 0-cycle on $X$, i.e., a sum of $d$ distinct
points, on a surface $X$.

{\bf Definition.} $\xi$ is said to be Caley-Bacharash (CB) with
respect to invertible sheaf $L$ on $X$ if for every point $p \in \xi$
the natural inclusion $0 \to H^0(X,J_{\xi}(L)) \to
H^0(X,J_{\xi-p}(L))$ is an isomorphism.

In other words, the pair $(\xi,L)$ is Caley-Bacharash if for every
point $p \in \xi$ and every curve $C$ (if any) in the linear system
$|L|$ on $X$ containing $\xi-p$ the point $p$ is also on $C$.

{\bf Lemma.} A Caley-Bacharash pair $(\xi,L)$ is special.

{\bf Proof.} Since $(\xi,L)$ is Caley-Bacharash, the exact sequence
$$
   \triple{J_{\xi}(L)}{J_{\xi-p}(L)}{L \tensor k(p)}
$$
induces the long exact sequence
$$
   0 \to L \tensor k(p) \to H^1(X,J_{\xi}(L)) \to H^1(X,J_{\xi-p}(L)) \to 0
$$
which implies $h^1(X,J_{\xi}(L)) = h^1(X,J_{\xi-p}(L)) + 1$
and therefore $\delta(\xi,L) = \delta(\xi-p,L) + 1 \ge 1$.

{\bf Examples.} Assume that the linear system $|L|$ is not empty.

1. If $\xi = [p], p \in X$, then $(\xi,L)$ is Caley-Bacharash iff $p$
is a base point for $|L|$.

2. If $\xi = [p] + [q]$ and $|L|$ is basepoint-free, then $(\xi,L)$ is
Caley-Bacharash iff $\phi_L(p) = \phi_L(q)$.

3. If $\deg \xi = 3$, $\xi = [p]+[q]+[r]$, and $|L|$ is basepoint-
free, then $(\xi,L)$ is Caley-Bacharash iff either $\phi_L(p) =
\phi_L(q) = \phi_L(r)$, or $\phi_L(p)$, $\phi_L(q)$ and $\phi_L(r)$
are {\it all distinct} and $\dim <p,q,r>_L = 1$. Note that there are
special degree 3 0-cycles which are not Caley-Bacharash.

4.  If $\deg \xi = d$ and $\xi$ is a sum of $d$ distinct points on $X$
which are not base points for $|L|$, then $(\xi,L)$ is Caley-Bacharash
iff $\phi_L(p) \in <\xi - p>_L $ for every $p \in \xi$.

\begin{for-me}
Full description:

either (a) $|L|$ is empty, or (b) $p$, $q$ and $r$ are base points for
$L$, or $\phi_L(p) = \phi_L(q) = \phi_L(r)$, or $\phi_L(p)$,
$\phi_L(q)$ and $\phi_L(r)$ are distinct and $\dim <p,q,r>_L = 1$,
where $<p,q,r>_L := <\phi_L(p),\phi_L(q),\phi_L(r)>$ is the linear
envelope of the points $\phi_L(p)$, $\phi_L(q)$ and $\phi_L(r)$ in the
projective space $|L|\dual$
\end{for-me}

\subsubsection{Speciality index of a $0$-cycle on a curve.}
 
Let $C$ be a (smooth) non-hyperelliptic curve and $\xi$ be a simple
0-cycle on $C$ of degree $d$, i.e., an effective divisor consisting of
$d$ distinct points.

{\bf Definition.} The speciality index of 0-cycle $\xi$ of degree $d$
with respect to the canonical linear system $|K_C|$ is defined by 
$$ \dim <\xi>_{K_C} = d - 1 - \delta_C(\xi,K_C)
$$  

{\bf Lemma.}
$$
   \delta_C(\xi,K_C) = \dim |\xi|
$$  

{\bf Proof.} This is equivalent to the Riemann-Roch theorem.

\subsubsection{Speciality index of a $0$-cycle on a curve on a surface.}
Let $X$ be a surface with $h^1(X,\ox) = 0$, $C$ be an effective
Cartier divisor on $X$ such that $\Omega^2_X \tensor_{\ox} \ox(C)$ is
ample, and let $\xi$ be a Cartier divisor on $C$. Let $L = \ox(C)$, $A
= \oc(\xi)$ and $B = N_{C/X} A^{-1}$.

{\bf Lemma.} 
\label{xi<C<X}
$$
    \delta_X(\xi,L) = h^1(C,B) - p_g(X)
$$

{\bf Proof.}  Let $i: C \to X$ be the canonical embedding. The
isomorphism $\alpha: \Ext^1_{\ox}(i_* B, \ox) \isom H^0(C,A)$ of
~\ref{sheafExt^1(B,O)=A} and the description of extensions given in
~\ref{0-1-is-torsion-free:global-case} give an extension 
$$
   \xyTriple{\ox}{\jk(L)}{i_*B}{s}{}
$$
with $s$ vanishing at $C$. The induced long exact
sequence
$$
   0 \to 
   H^1(X, \jk(L)) \to
   H^1(X, i_* B) \to
   H^2(X, \ox) \to
   H^2(X, \jk(L))
   \to 0
$$
gives 
$ h^1(\jk(L)) = h^1(X,B) - p_g + h^2 \jk(L) =
                h^1(C,B) - p_g + h^2 (L) =
                h^1(C,B) - p_g 
$
,where we used the Kodaira vanishing theorem.

\subsubsection{$0$-cycle on a curve on $K3$}

Let $X$ be a K3 surface, $C$ be a smooth irreducible curve on $X$ such
that $L = \ox(C)$ is ample, and let $\xi$ be a Cartier divisor on
$C$. Let $A = \oc(\xi)$ and $B = K_C A^{-1}$.

{\bf Lemma.} 
$$
    \delta_X(\xi,L) = h^0(C,A) - 1
$$

This is a corollary of the previous lemma. % ~\ref{xi<C<X}

If $C$ is non-hyperelliptic and $\xi$ is a simple 0-cycle, then the
lemma above can be reformulated as

$$
    \delta_X(\xi,L) = \delta_C(\xi,K_C)
$$

We also give another proof:

Let $s$ be the canonical section of $L=\ox(C)$. If $\phi_L$ is a
morphism associated with linear system $|L|$ and $H$ is a hyperplane
in $H^0(X,L)\dual$ orthogonal to $s$, there is a diagram
$$
   \xymatrix
   {
      C \ar[d]^{i}\ar[r]^-{(\phi_L)|_C} & H \ar[d] \\
      X           \ar[r]^-{\phi_L}      & \PP H^0(X,L)\dual
   }
$$
Consider the short exact sequence
$$
   \xyTriple{\ox}{L}{L|_C}{s}{}
$$
and the associated exact sequences
$$
   0 \to
   H^0(X,\ox) \strelka{s} 
   H^0(X,L) \to
   H^0(C,L|_C) \to 0
$$
and
$$
   0 \to
   H^0(C,L|_C)\dual \to
   H^0(X,L)\dual \to
   H^0(X,\ox)
   \to 0
$$
It is clear that there is an isomorphism $\PP H^0(C,L|_C) \dual \to H$
which makes the following diagram commutative:
$$
   \xymatrix
   {
      C \ar[r]^-{\phi_{(L_C)}} 
        \ar[rd]_-{(\phi_L)|_C} & 
      { \PP H^0(C,L|_C) \dual } \xyIsom[d] \\
      {} & H
   }
$$

Since the restriction $L|_C \simeq \Omega^1_C$, the 
linear spans $<\xi>_{K_C}$ and $<\xi>_{L}$ coincide, and
the lemma follows.
 
          % no extra [], spell-checked
          % \subsection{(*) Special 0-cycles on surfaces}
          % ----------

% -------------------------------- gg things ---------------------------

\section{Injectivity of generic evaluation map for globally
               generated vector bundles}

\providecommand{\trivial}[1]{#1 \tensor \OO_X}

\begin{for-me}
\subsection{Three trivial examples.}

1a: If $E$ is a locally free sheaf, $V \subset H^0(E)$
and $\phi=\phi_V: V \tensor \OO_X \to E$ is a canonical homomorphism,
then $\phi$ may be non-injective. Example:
$$\trivial{H^0(\mathbb P_1,\OO(1))} \to \OO(1) \to 0$$

1b: If $E$ is a locally free sheaf, $v= \dim V$, and $v \le \rk E$,
then it is possible that for all $V \in Gr(v,H^0(E))$ the map $\phi_V$
is not injective. For example, let $E=\OO_{\PP_1}(-1)+\OO_{\PP_1}(5)$,
and $\dim V = 2$. 

1c: If $E$ is locally free and globally generated, $v \le \rk E$, then
it is possible that $\phi_V$ is not injective. Example: let $E$ be a
globally generated locally-free sheaf, $L \subset E$ be a locally free
subsheaf of rank one, and let $V \subset H^0(L) \subset H^0(E)$, $\rk
V > 1$. Then $\phi_V:  \trivial{V} \to E$ is not injective.
\end{for-me}

\subsection{The injectivity lemma.}
\begin{thing}
\label{injectivity-lemma}
{\bf Lemma.} Let $X$ be a reduced and irreducible scheme, and $E$ be a
locally free sheaf on $X$ generated by global sections. Let $v \le \rk
E$. Then there is a {\it nonempty} open subset in the Grassmanian variety
$Gr(v,H^0(E))$ such that for all $V$ in this subset the map $\phi_V$
is injective.
\end{thing}

{\bf Proof.} Let $H=H^0(E)$. Consider the diagram
\begin{equation}
   \xymatrix@C=7mm@R=5mm
   {  
    0\ar[r] & {\ker\psi} \ar[r] & \trivial{H} \ar[r]^{\psi} & E \ar[r] & 0 \\
    0\ar[r] & {\ker \phi_V}\ar[r] \xyArrow{\cup}{u} & \trivial{V} \ar[r]^{\phi_V}            \xyArrow{\cup}{u} & E \xyRavno[u] & {}
    }
   \label{diagram:E-H-V}
\end{equation}

  Note that $\ker \psi$ is locally free and $\ker \phi_V = (\ker \psi)
\intersect (V \tensor \OO_X)$.

Let $K$ be the field of rational functions on $X$. For a sheaf $F$ on
$X$ we denote by $F(K)$ the $K$-vector space of rational sections of
$F$, which is isomorphic to the fiber of $F$ at the generic point of
$X$. The diagram above induces the diagram
$$
   \xymatrix@C=7mm@R=5mm
   { 
     0 \ar[r] & (\ker \psi)(K)  \ar[r] & H \tensor_k K \ar[r] & E(K)\ar[r] & 0 \\
     0 \ar[r] & (\ker \phi_V)(K)\ar[r] \xyArrow{\cup}{u} & V\tensor_k K                      \xyArrow{\cup}{u}  & {} &
   }
$$
Note that $(\ker \phi_V)(K)$ is the intersection of $K$-vector subspaces
$(\ker \psi)(K)$ and $V \tensor K$ in $H \tensor K$.

For a $k$-vector space $V$ let $V_K = V \tensor_k K$. Let $A =(\ker
\psi)(K)$. 
% Let $\dim_k H = u$, $\dim_k V = v$ and $\rk E=r$.  
The inequality $\dim V \le \rk E$ implies that
$$
   \dim_K A + \dim_K V_K \le \dim_K H(K)
$$

Now the ``Grassmanian Lemma'' ~\ref{grassmanian-lemma} below implies
that for a generic $V \subset H$ the intersection $V_K \cap A = 0$,
i.e., $(\ker \phi_V)(K) = 0$. But since $\ker \phi_V$ is a subsheaf of
a torsion-free sheaf $V \tensor \OO_X$, we have $\ker \phi_V = 0$.

\subsection{Grassmanian lemma}

\providecommand{\bG}{\mathbf{G}} % bold G
\providecommand{\bGr}{\mathbf{Gr}} % bold Gr
\providecommand{\bX}{\mathbf{X}} % bold X
\providecommand{\bSigma}{\mathbf{\Sigma}} % bold \Sigma
\providecommand{\bH}{\mathbf{H}} % bold H
\providecommand{\bI}{\mathbf{I}} % bold I
\providecommand{\bL}{\mathbf{L}} % bold I
\providecommand{\bOmega}{\mathbf{\Omega}} % bold \Omega
\begin{thing}
\label{grassmanian-lemma}
{\bf Lemma.}
Let $K/k$ be an extension of fields, $F$ be a vector
space over $k$, and let $\bG = \bGr_k(v,F)$ be the Grassmanian variety of
$v$-dimensional $k$-vector subspaces in $F$. Let $F_K = F \tensor_{k}
K$ and let $A \subset F_K$ be a $K$-vector subspace of dimension $a$
such that $a + v \le \dim_k F$. Let $X_A$ be the set of all
$v$-dimensional $k$-subspaces $V$ in $F$ such that $\dim_K (V_K
\intersect A) > 0$. Then there is a proper algebraic subvariety
$\bX_A$ in $\bG$ such that $X_A$ is the set of $k$ - points of
$\bX_A$.
\end{thing}

 {\bf Proof.}  Let $\Gr_k(v,F)$ be the set of all $v$-dimensional
$k$-vector subspaces in $F$; we have $\Gr_k(v,F) = \bG(k)$.  Let
$\Gr_K(v,F_K)$ be the set of all $v$-dimensional $K$-vector subspaces
in $F_K$ and the $\bGr_K(v,F_K)$ be the corresponding Grassmanian
variety. There is a canonical isomorphism $\bGr_K(v,F_K) \simeq \bG
\tensor_k K$. (Note that the results of ~\cite{ega1}, pp9.7-9.8 can be
applied to the Grassmanian of subspaces of $F$ since the base scheme
$S = \Spec k$).

Let $t = \dim_k F$ and $m = t-a-v+1$; we have $m \ge 1$.  Let
$$
   \Sigma_A = \{W \subset F_K: \dim_K W =v, \dim_K(W \intersect A) >
   0 \}
$$ It is well-known that $\Sigma_A$ is the set of $K$-points of a
proper algebraic subvariety $\bSigma_A$ in $\bG_K$ of codimension
$m$. (The class of $\bSigma_A$ in the cohomology ring or Chow ring is
a special Schubert class $\sigma_m = \{m,0,\dots,0\}$,
cf. ~\cite{fulton}, p. 14.7).

  Let $\PP = \PP_k(\Lambda^vH)$ and $\PP_K = \PP \tensor_k K$; there
is a natural isomorphism $\PP_K \simeq \PP_K(\Lambda^v H_K)$.  There
is a commutative diagram
$$
  \xymatrix
  {
    \bGr_K(v,H_K) \ar[d]\ar[r]^{p_K} & \PP_K(\Lambda^v H_K)\ar[d] \\
    \bGr_k(v,H)   \ar[r]^{p}         & \PP_k(\Lambda^v H) 
  }
$$
where $p$ and $p_K$ are Plucker morphisms.

   Consider first the case $t = a + v$, $m=1$. In this case $\bSigma_A
= \bG_K \intersectIn{\PP_K} \bH_A$ for some hyperplane section $\bH_A
\subset \PP_K$. 

  Let $H_A = \bH_A(K)$, $G = \bG(k)$, $G_K = \bG_K(K)$, $P = \PP(k)$
and $P_K = \PP_K(K)$. There are natural embeddings of $G$ into $G_K$
and $P$ into $P_K$. Note that the intersection $H_A \intersectIn{P_K} P$
is the set of $k$- points of some linear subspace $\bI \subset \PP$.
 We have
\begin{multline*}
   X_A = 
   \Sigma_A  \intersectIn{G_K} G
   = 
   (H_A \intersectIn{P_K} G_K) \intersectIn{G_K} G 
   \stackrel{[1]}{=}
   (H_A \intersectIn{P_K} P) \intersectIn{P} G
   = 
   I \intersectIn{P} G
   = \\ =
   \bI(k) \intersectIn{\PP(k)} \bG(k) 
   =
   (\bI \intersectIn{\PP} \bG) (k)
\end{multline*}
where the equality [1] is induced by the diagram 
$$
   \smallXyMatrix
   {
     P_K \xyArrow{\supset}{r} &  G_K \\
     P   \xyArrow{\supset}{r}\xyArrow{\cup}{u} &  G \xyArrow{\cup}{u}
   }
$$
Therefore we can let $\bX = \bI \intersect \bG$.

Note that if $A \subset F_K$ is in a general position and $\dim K/k >
\dim_k \bG$, then $X_A$ is empty.

In the general case consider the canonical map
$$
   u_A:   A \tensor_k \Lambda^{v-1} F \to \Lambda^v F
$$
and let $\bL_A = \PP(\Im(u_A))$. $\bL_A$ is a linear subspace in $\PP$,
and there is a fibered square
$$
   \smallXyMatrix
   {
     \bG        \ar[r]                  &  {\PP} \\
     \bSigma_A  \ar[r]\xyArrow{\cup}{u} &  {\bL_A} \xyArrow{\cup}{u}
   }
$$
which allows one to repeat the computation above.

{\bf Example.} Let $\dim_k F = 4$, $\bG =\bGr(2,F)=\bGr(1,\PP(F))$,
and let $A$ be one-dimensional vector subspace in $F$. We have $L_A
\subset \PP(\Lambda^2 F)$, $L_A \simeq \PP(A \tensor_k (F/A))$, 
$\dim L_A = 2$. Note that $L_A \subset \bG$, and therefore $\Sigma_A =
L_A \intersect G = L_A$. If $p = \PP(A)$, $p \in \PP(F)$, then
$\bSigma_A = \bSigma_p$ is the variety of lines in $\PP(F) \simeq
\PP^3$ containing the point $p$.

   Let $L_A = H_1 \intersect H_2 \intersect H_3$, where $H_1$, $H_2$
and $H_3$ are hyperplanes in $\PP(\Lambda^2 F)$. Then the intersection
of $H_1 \intersect H_2$ and $\bG$ in $\PP(\Lambda^2(F))$ is proper,
and $G \intersect H_1 \intersect H_2$ has two irreducible components,
one of them $G \intersect H_1 \intersect H_2 \intersect H_3$. If $l_1$
and $l_2$ are two lines in $\PP(F)$ intersecting at the (single) point
$p$, and if $H_1 = \bSigma_{l_1}$ and $H_2 = \bSigma_{l_2}$, then
$$\bSigma_{l_1} \intersect \bSigma_{l_2} = \bSigma_p + \bOmega_e,$$
where $\bOmega_e$ is the variety of lines contained in the plane $e$
spanned by $l_1$ and $l_2$. This decomposition corresponds to the
equality $\sigma_1^2 = \sigma_2 + \sigma_e$ in the Chow ring of $\bG$,
(see, for example, ~\cite{fulton}, Example 14.7.2). In coordinates, we have
$\bG = a_{12}a_{34}-a_{13}a_{24}+a_{14}a_{23}=0$, and if $H_1=
(a_{12}=0)$, $H_2= (a_{13}=0)$ and $H_3= (a_{14}=0)$, then 
$\bG \intersect H_1 \intersect H_2 = (a_{12}=0,a_{13}=0,a_{14}=0) +
(a_{12}=0,a_{13}=0,a_{23}=0) = G \intersect H_1 \intersect H_2
\intersect H_3 + \Omega_e$.

%\end{document}

%           no extra [], spell-checked
%           \section{(.) Injectivity of generic evaluation map for globally
%                      generated vector bundles}
%           \subsection{The injectivity lemma.}
%           \subsection{Grassmanian lemma}
%           -----------------------------

\newcommand{\Proj}{\mathop{\rm Proj}}
\newcommand{\Sing}{\mathop{\rm Sing}}
\newcommand{\Sym}{\mathop{\rm Sym}}
\newcommand{\lfE}{{\cal E}} 

\section{One application of Bertini theorem}
\subsection{Smoothness lemma}
\label{lemma-on-cotangent-sheaf}

I did not find a reference for the following lemma, though it is
probably standard:

\begin{thing}
{\bf Lemma.} Let $Y$ be a smooth irreducible affine scheme over a
field $k$, $R = \Gamma(Y,\OO_Y)$, $s = (s_1,\dots,s_r)$ be a sequence
of elements in $R$, and let $Z=Z(s)$ be a subscheme in $Y$ given by
$s$.

(a) There is a complex of coherent sheaves on $Z$
$$
   \OO_Z^r 
   \strelka{(ds_1|_Z, \dots, ds_r|_Z)}
   (\Omega^1_{Y/k})|_Z 
   \to
   \Omega^1_{Z/k} 
   \to 0
$$

(b) If
$ds_1|_Z, \dots, ds_r|_Z$ are $\OO_Z$-independent sections of the
locally free sheaf $(\Omega^1_{Y/k})|_Z = \Omega^1_{Y/k} \tensor \OO_Z$,
(i.e., if $(ds_1 \wedge \dots \wedge ds_r) \tensor 1$ is a
non-vanishing section of $\Omega^r_{Y/k} \tensor \OO_Z$),
then $(s_1,\dots,s_r)$ form a regular sequence in $R$, and
$Z$ is smooth over $k$ of dimension $\dim _k Y - r$.
\end{thing}

{\bf Proof.} Let us prove (b). Consider first the case $r = 1$. In this case $Z$ is
given by a single element $s \in R$. Let $j: Z \to X$ be the canonical
embedding, and let $J$ be the ideal generated by $s$ in $R$.  Taking
the global sections of the exact sequence
$$
   N_{Z/X} \to j^*(\Omega^1_{X/k}) \to \Omega^1_{Z/k} \to 0
$$ 
(cf. ~\cite{ega4}, 4.16.4.21), we get the exact sequence of $R/J$-modules
$$
   J/J^2 \strelka{d} \Omega^1_{R/k} \tensor_R R/J \to \Omega^1_{(R/J)/k} \to 0
$$
(cf. ~\cite{ega1}, 0.20.5.12.1), where $J/J^2$ is a free $R/J$-module of rank
one generated by $s \mod J^2$.  Let $p \in Z$. It is clear that the
condition $\Tor_1^{R/J}(\Omega^1_{(R/J)/k}, k(p)) = 0$ is equivalent
to $ds \tensor k(p) \ne 0$, which follows from the assumption of the lemma.  It
follows that $\Omega^1_Z$ is a locally free sheaf of rank $d-1$, where
$d = \dim Y$, and therefore $Z$ is nonsingular of dimension $d-1$.

It is clear that $\{s\}$ forms a regular sequence: indeed, if $st=0$
for some $t \in R$, we have $sdt + tds=0$, and $X = V(s) \cup
V(t)$. Therefore we have $tds=0$ on $V(s)=Z$, which contradicts to the
assumption on $ds|Z$. (We do not want to use the irreducibility of $X$
in order to be able to use this argument in the induction step later.)

 Consider now the case $r > 1$.  Let $R_i = R/(s_1,\dots,s_i)$, $i =
1,\dots,r$; we have $R_i = R_{i-1}/(s_i)$. Let $Z_i = \Spec R_i$.

Consider the map $i_r: \OO_Y^r \to \Omega^1_{Y/k}$ given by
$(ds_1,\dots,ds_r)$. Taking the $r$-th exterior power, we get
$$
   \wedge^r i_r: \wedge^r \OO_Y^r \to \Omega^r_{Y/k}
$$ Since $(ds_1 \wedge \dots \wedge ds_r)|_Z \ne 0$, $\wedge^r i_r$
gives a section of $\Omega^r_{Y/k}$ not vanishing on $Z$, and
therefore there is an affine subscheme $U \subset Y$ such that $Z
\subset U$ and $ds_1,\dots, ds_r$ are $\OO_U$-independent sections of
$\Omega^1_Y|_{U} \isom \Omega^1_U$. Let us substitute $Y$ by $U$.

Assume that for some $i\ge 2$ $Z_{i-1}$ is smooth of dimension
$d-(i-1)$ and that the sequence $(s_1,\dots,s_{i-1})$ is regular.
Now $Z_i$ is given by a single equation $s_i$ in $R_{i-1}$,
and the exact sequence of $R_i$-modules
$$ 
  (R_{i-1} s_i)/(R_{i-1} s_i^2)
  \strelka{d} 
  \Omega^1_{R_{i-1}/k}
  \tensor_{R_{i-1}} R_i \to \Omega^1_{R_i/k} \to 0
$$
gives that $\Omega^1_{R_i/k}$ is locally free at $p \in Z_i$ of rank
$d -i$ iff $ds_i \tensor k(p) \ne 0$. But since 
$(ds_1 \wedge \dots \wedge ds_r) \tensor k(p) \ne 0$ for all $p \in Y$, 
the lemma follows.

 We can repeat the argument for $r=1$ to see that $s_i$ is not a zero
divisor in $R_{i-1}$.

  The proof of (a) is straightforward.

%--------------------------------------------------------------------------------

\subsection{$Z(s)$ for a generic section of a globally generated
bundle}

For a section $s$ of a locally free sheaf $\lfE$ we denote by $Z(s)$
the scheme of zeroes of $s$.

\begin{thing}
\label{generic-section-of-a-gg-lf-sheaf}
{\bf Proposition.}  Let $X$ be a smooth proper scheme over a field $k$
of characteristics 0, and let $\lfE$ be a locally free sheaf on $X$
generated by global sections. Then

   (a) If $\rk E > \dim X$, then a generic sections $s \in H^0(E)$ does
not vanish on $X$; 

   (b) If $\rk E \le \dim X$, then for a generic $s \in H^0(E)$ the
scheme $Z(s)$ is smooth of codimension $r$ in $X$, and, moreover, each
point $q \in Z(s)$ has an affine neighborhood on which $s$ is given
by regular sequence.  
\end{thing}

{\bf Proof.} Let $E$ be the vector bundle $\Spec \mathop{\rm Sym} \lfE
\dual$.  Consider the diagram
$$
   \xymatrix
   {
     {\PP(E\dual)} \ar[d]^{\pi} \ar[rd] & {} \\
     X \ar[r] & {\Spec k}
   }
$$ where $\PP(E\dual)$ is the notation of ~\cite{fulton} for the
projectivization of a vector bundle (in the notation of EGA2 
this is $\PP(\lfE) = \Proj \mathop{\rm Sym} \lfE$).

Let $\OO_{E\dual}(1)$ be the canonical invertible sheaf on
$\PP(E\dual)$. Since $\pi_* \OO_{E\dual}(1) \simeq \lfE$, the diagram
above induced the isomorphism
$$
    A: H^0(X,\lfE) \isom H^0(\PP(E\dual), \OO_{E\dual}(1))
$$
If $s \in H^0(X,\lfE)$, let $D_s = Z(A(s))$ be the divisor of zeroes
of the corresponding section of $\OO_{\PP(E\dual)}(1)$.

\begin{thing}
{\bf Lemma.} The locally free sheaf $\lfE$ is generated by global
sections iff the linear system $|\OO_{E\dual}(1)|$ on $\PP(E\dual)$ is
basepoint-free.
\end{thing}

{\bf Proof.} Assume that the linear system $ |\OO_{E\dual}(1)|$ has a
base point $p \in \PP(E\dual)$, i.e., that $p \in D_s$ for all for all
$s \in H^0(X,E)$. Let $q = \pi(p)$, and let $H_p \subset E(q)$ be the
hyperplane in the fiber of the vector bundle $E$ over the point $q$
orthogonal which consists of vectors orthogonal to the point $p \in
E\dual(q)$. The condition $p \in D_s$ is equivalent to the condition
$s(q) \in H_p$. Since this should be true for all sections $s \in
H^0(X,E)$, it follows that the fiber $E(q)$ is not generated by the
global sections of $E$.

The inverse statement of lemma is
straightforward. 

\begin{thing}
{\bf Lemma.} The divisor $D_s$ is singular iff $Z(s)$ is neither smooth
of dimension $d-r$ nor empty.
\end{thing}

{\bf Proof.}
Let $q \in X$, and let $U$ be an affine subset in $X$ containing $q$
such that both $\lfE$ and $\Omega^1_X$ can be trivialized when
restricted to $U$. Let us choose an isomorphism $\alpha_U: \lfE|_U \to
\OO_U^r$, and let $(e_1,\dots,e_r)$, $e_i \in \Gamma(U,\lfE|_U)$, be
the preimages of the canonical sections of $\OO_U^r$.

Taking the dual of $\alpha_U$, we get a trivialization $\alpha_U\dual:
\OO_U^r \isom \lfE\dual|_U$ and, therefore, an isomorphism $U \times
\PP^{r-1} \isom \PP(E\dual|_U)$. Let $(e_1\dual,\dots, e_r\dual)$ be
the dual trivialization of $\lfE\dual|_U$, and let $x_1,\dots,x_r$ be
the canonical homogeneous coordinates on $\PP^{r-1}$.

  Let $R = \Gamma(U,\OO_U)$ and let us make a base change to $U$. We
have $s = \sum_1^r s_i e_i$ for some $s_i \in R$. It is clear that
$D_s \subset \PP(\lfE\dual)$ can be given by by the equation
$$
   f_s = \sum_1^r s_i x_i
$$

The singularities of $D_s$ are given by the condition $df_s = \sum_1^r
s_i dx_i + \sum_1^r x_i ds_i = 0$, which is equivalent to the
condition (1) $s_i=0$ for all $i$ and (2) $\sum_1^r x_i ds_i = 0$.

  In other words, let $Z = Z(s)$, and consider the map of the locally
free sheaves $\phi: \OO_X^r \to \Omega^1_X$, where $\phi =
(ds_1,\dots,ds_r)$. Consider the restriction of $\phi$ on $Z$ and
consider the corresponding map of vector bundles
$$ 1_Z^r \strelka{\phi|_Z} T_X\dual |_Z
$$ 
Let $C_s$ be the preimage of the zero section of $T_X\dual |_Z$ under
$\phi_Z$, and let $PC_s$ be the projectivization of $C_s$.  We have
$PC_s = \Sing (D_s)$.

(Scheme-theoretically, for a morphism of locally free sheaves ${\cal E}
\to {\cal F}$ the preimage of the zero section of $F$ can be defined
as $ \Spec ( \left.  \Sym {\cal E}\dual \right/ \phi( \Sym^{+} {\cal
F}\dual ) \cdot \Sym {\cal E}\dual ) $ and the corresponding
projective cone as a $\Proj$ of the same graded sheaf of algebras.)

Now lemma ~\ref{lemma-on-cotangent-sheaf} gives a complex
$$
   \OO_Z^r 
   \strelka{(ds_1|_Z, \dots, ds_r|_Z)}
   \Omega^1_{X}|_Z 
   \to
   \Omega^1_{Z} 
   \to 0
$$

  Let now $q \in x$, and consider the exact sequence
$$ 
   0 \to
   \ker \phi \tensor k(q) 
   \to
   k(q)^r
   \strelka{\phi \tensor k(q)}
   \Omega^1_X \tensor k(q) 
   \to
   \Omega^1_Z \tensor k(q) 
   \to 0
$$
% It is clear that $\dim_k \ker \phi \tensor k(q) > 0$ iff 
$\dim_k \Omega^1_Z \tensor k(q) > r-d$. If $C_s(q)$ is a
scheme-theoretic fiber of $C_s$ over $q$, then set of $k$-points is in
a natural correspondence with $\ker \phi \tensor k(q)$.

Let us come back to the general case of not-necessarily affine $X$.
Let $\Sing_{d-r}(Z)$ be the set of points $q \in Z$ for which $\dim_k
\Omega^1_Z \tensor k(q) > d-r$, i.e., the set of points at which $Z$
fails to be smooth of expected codimension. It is clear that this is a
set of points of a closed subscheme in $Z$, and that the canonical
morphism $PC_s \to X$ factorizes through $\Sing_{d-r}(Z)$. In another
words, there is a diagram
$$
   \xymatrix
   {
      {\PP(E\dual)} \ar[d]\xyArrow{\supset}{r} & D_s \xyArrow{\supset}{r} & {\Sing D_s = PC_s} \ar[d]^{\pi'} \\
      X \xyArrow{\supset}{r} & Z  \xyArrow{\supset}{r} & \Sing_{d-r}(Z) 
   }
$$ 
It is clear that $\pi'$ is epimorphic, and the statement of the lemma
follows.

% lemma 3
\begin{thing}
{\bf Lemma.}
If $r > d$, then either $Z(s)$ is empty and $D_s$ is smooth, or $Z(s)$
is not empty and $D_s$ is singular. If $r \le d$, then either $Z(s)$
is empty or smooth of expected dimension $d-r$ while $D_s$ is smooth,
or $Z(s)$ fails to be smooth of expected dimension while $D_s$ is
singular.
\end{thing}

This lemma follows from the previous one immediately.

Let us now use the assumption $\mathop{\rm char} k = 0$ and use the
Bertini's theorem (cf. ~\cite{hartshorne}):

% Corollary 1
\begin{thing}
{\bf Corollary.} If $r > d$, then the Bertini theorem implies that
$Z(s)$ is empty for the generic $s \in H^0(X,\lfE)$, and if $r \le d$,
then the Bertini theorem implies that $Z(s)$ is either empty or smooth
of dimension $d-r$ for the generic $s \in H^0(X,\lfE)$, and the
generic section $s$ gives a regular sequence in the coordinate ring of
some affine neighborhood of every point $q \in X$. 
\end{thing}

This finishes the proof of the theorem.

%\end{document}

%           no extra [], spell-checked
%           \section{Some applications of the Bertini theorem}
%              \subsection{Smoothness lemma}
%              \subsection{$Z(s)$ for a generic section of  g. generated bundle}
%           -------------   

\section{Computation of $\jk(L)/ \OO_X$ }

\providecommand{\oc}{\OO_C}
\providecommand{\ox}{\OO_X}

Let $X$ be a surface, $\xi$ be a 0-subscheme in $X$, $L$ be an
invertible sheaf on $X$ and $s \in H^0(X,L)$ be $\ox-regular$
section of $L$. Considering the embedding $\jk(L) \subset L$, we get
an $\ox$-regular global section of $L$. Let $C$ be the Cartier
divisor of zeroes of this section and let $i: C \to X$ be
the embedding of the support of $C$ into $X$. Twisting the homomorphism
\begin{equation}
  s: \ox \to \jk(L)
  \label{section-s}
\end{equation}
with $\ox(-C) \simeq L^{-1}$, we see that $\xi$ is a Cartier divisor
on $C$. Let $l = i^*(L)$.

\begin{thing}
\label{0-1-lemma}
{\bf Lemma.} The map ~\ref{section-s} can be included into the exact
sequence
$$
   \triple{\ox}{\jk(L)}{i_*(l(-\xi))}
$$
\end{thing}

{\bf Proof.} Let $j: \xi \to C$ be the canonical embedding of the support
of $\xi$ to $C$ and let $k = ij$. Then the inclusions
$$
   \xi \stackrel{j}{\subset} C \stackrel{i}{\subset} X
$$
give the filtration $J_{C,X} \subset J_{\xi,X} \subset \ox$ which induces
the diagram
$$
 \displayF{J_{C,X}}{J_{\xi,X}}{\ox}{k_*\ok}{i_*J_{\xi,C}}{i_*\oc}{i_*j_*\ok}
$$
This diagram can be rewritten as
$$
  \displayF{\ox(-C)}{J_{\xi,X}}{\ox}{k_*\ok}{i_*\oc(-\xi)}{i_*\oc}{i_*j_*\ok}
$$
twisting it with $\ox(C) \simeq L$, we get the diagram
$$
  \displayF{\ox}{\jk(L)}{L}{(k_*\ok)(L)}{i_*(l(-\xi))}{i_*l}{i_*((j_*\ok)(l))}
$$
left column of which gives the desired exact sequence.

%    no extra [], spell-checked
%    \section{Computation of $\jk(L)/ \OO_X$ }
%    -------------

\section{Good degeneration of $e_V$ in the globally generated case}
\begin{thing}
{\bf Proposition.}
\label{good-defeneration-in-the-gg-case}
  Let $X$ be a smooth surface over a field $k$ of characteristics 0,
and $E$ be a globally generated locally free sheaf on $X$.  Let $l$ be
an integer, $1 \le l \le \rk E$. Then for the generic vector subspace
$V \subset H^0(X,E)$ of dimension $l$ the cokernel of the evaluation
map $$e_V: H^0(X,E) \tensor_{k} \OO_X \to E$$ is acceptable. 
\end{thing}

{\bf Proof.}

 By Lemma ~\ref{injectivity-lemma}, for a generic $V \subset H^0(E)$
of dimension $l$ the map $e_V$ is injective. Let $V \subset H^0(E)$ be
generic, and let $l = \dim_k V$.

(a) Assume first that $l \le \rk E - 2$. In this case $r - l \ge
2$. By lemma ~\ref{generic-section-of-a-gg-lf-sheaf}, for generic $v
\in V$ the corresponding section $s_v \in H^0(X,E)$ does not vanish on
$X$, and therefore the factor-sheaf $E' = E/ e_V(v \tensor \OO_X)$ is
locally free. There is a diagram 
$$
   \displayE{v \tensor \OO_X}{V \tensor \OO_X}{E}{F}{(V/v) \tensor
\OO_X}{E'} 
$$ 
Applying induction to the pair $(E',V/v)$, we conclude that $F$ is
locally free of rank $r-l$.

(b) Assume now that $\dim_k V = \rk E - 1$. Following the step (a), we
can choose a generic $l-1$ -dimensional vector subspace $W \subset V$
such that $E' = E/ e_V(W \tensor O)$ is a locally free sheaf of rank 2
which can be included into the diagram
$$
   \displayE{W \tensor \OO_X}{V \tensor \OO_X}{E}{F}{(V/V) \tensor
             \OO_X}{E'}
$$
 Consider the pair $(E',V/W)$. $E'$ is a locally free sheaf of rank 2
and $V/W$ is one-dimensional vector subspace in $H^0(X,E')$. Let $u
\in V/W$, $u \ne 0$, and let $s_u$ be the corresponding section of
$E'$. The proposition ~\ref{generic-section-of-a-gg-lf-sheaf} ensures
that by choosing $W$ to be generic enough we have that either $Z(s)$
is empty or that each point $x \in X$ has an affine neighborhood $U
\simeq \Spec R$ such that $s|_U$ is given by a regular sequence
$(s_1,s_2)$, $s_i \in R$.  Let $J = Rs_1 + Rs_2$; the (cohomological)
Koszul resolution
$$
   0 \to 
   R \strelka{(s_1,s_2)}
   R^2 \strelka{(-s_2,s_1)}
   J
   \to 0
$$
proves that $F|_U \simeq J$, which finishes the proof. 

    (c) Assume now that $\dim V = \rk E$. By choosing a generic $l-1$
-dimensional vector subspace $W \subset V$ and considering $E' = E/ (W
\tensor \OO_X)$ as in (b) we get a torsion-free sheaf $E'$ of rank 1
and 1-dimensional subspace in the vector space $H^0(X,E')$. Since $E'$
is torsion free, it is of the form $\jk(L)$ for a 0-subscheme $\xi
\subset X$, and by lemma ~\ref{0-1-lemma}, $F$ in the diagram above is
isomorphic to the direct image of an invertible sheaf on some Cartier
divisor $C$ on $X$, which completes the proof.

%      no extra [], spell-checked
%      \section{Good degeneration of $e_V$ in the globally generated case}
%      ----------

\section{ Construction of Brill-Noether loci}
\newcommand{\E}{{\cal E}} 
\newcommand{\F}{{\cal F}} 
\newcommand{\G}{{\cal G}} 
\subsection{Complex $K(E)$}
   $E$ be a torsion-free coherent sheaf on a surface $X$.  We
construct a complex of vector spaces which computes $H^i(X,E)$.  Let
$C$ be an effective divisor $X$ which is ample enough so that
$H^i(X,E(C))=0$ for $i = 1,2$. Let $i_1: C \to X$ be the embedding of
$C$ to $C$. Tensoring the adjunction sequence
$$ 
   \triple{\ox}{\ox(C)}{(i_1)_* \oc(C)}
$$
with $E$, we get 
$$
   \triple{E}{E(C)}{(i_1)_* E(C)|_C}
$$

Let $F= E(C)|_C$. We can choose an effective divisor $\xi$ on $C$
such that $H^1(C, F(\xi)) = 0$. Let $i_2: \xi \to
C$ be the embedding. We get a complex
$$
    \triple{F}{F(\xi)}{(i_2)_* F(\xi)|_{\xi}}
$$

Thus we get exact sequences 
\begin{gather*}
   0 \to H^0(X,E) \to H^0(X,E(C)) \strelka{\alpha}
          H^0(C,E(C)|_C) \to H^1(X,E) \to 0, \\
   0 \to H^1(C,E(C)|_C) \to H^2(X,E) \to  0, \\
   0 \to H^0(C,F) \strelka{\beta} H^0(C,F(\xi)) \strelka{\gamma} 
         H^0(\xi,F(\xi)|_{\xi}) \to H^1(C,F) \to 0
\end{gather*}

Let $K(E)$ be the complex of vector spaces 
$$
  K(E) = \left( 
                  0 \to H^0(X,E(C)) \strelka{\beta \circ \alpha}
                  H^0(C,F(\xi)) \strelka{\gamma}
                  H^0(\xi,F(\xi)|_{\xi}) \to 0
        \right)  
$$
which computes $H^*(X,E)$.

Note that this inductive procedure of choosing a flag of divisors
works on a variety of any dimension $d$ and gives a complex of length
$d$.

\subsection{Relative case}
Let $\E$ be a flat family of sheaves with base $S$ on a scheme $X$.
Consider the diagram
$$
 \xymatrix
 {
   X \times S  \ar[d]^{f_2}        & 
   C \times S  \ar[ld]_{f_1} \ar[l]_{j_1}      & 
   {\xi \times S}   \ar[dll]^{f_0}    \ar[l]_{j_2}   \\
   S & {} &  {} 
 }
$$

Let $p_X: X \times S \to X$ and $p_C: C \times S \to X$ be the
projection morphisms. We denote the lift of $C$ to $X \times S$ as
$C$, etc.  Let $\F = \E(C)|_{S \times C}$.  We can choose $C$ and
$\xi$ such that $R^n{f_2}_* \E(C) = R^n{f_1}_* \F(\xi) = 0$ for $n >
0$. Let $\G = \F(\xi)|_{\xi \times S}$.

  The exact sequences
$$
   \triple{\E}{\E(C)}{\E(C)|_{C \times S}}
$$
and 
$$
    \triple{ \F }{ \F(\xi) }{ \F(\xi)|_{\xi} }
$$
give long exact sequences
\begin{gather*}
    0 \to 
    {f_2}_* \E       \to
    {f_2}_* \E(C)   \strelka{\alpha}
    {f_2}_* {i_1}_* \E(C)|_{C \times S} \to
    R^1{f_2}_* \E  \to 0,  \\
    0 \to 
    R^1{f_2}_* {j_1}_* \F \to
    R^2{f_2}_* \E    \to  0,    \\
    0 \to 
    {f_1}_* \F        \strelka{\beta} 
    {f_1}_* \F(\xi)   \strelka{\gamma} 
    {f_1}_* {j_2}_* F(\xi)|_{\xi \times S} \to
    R^1{f_1}_* F    \to 0 
\end{gather*}
and a complex of $\OO_S$-modules
$$ 
   K(\E) =
                 \left( 
                     0 \to 
                     {f_2}_* \E(C)
                     \strelka{\beta \circ \alpha}
                     {f_1}_* \F(\xi) 
                     \strelka{\gamma}
                     {f_0}_* \G
                     \to 0
                 \right) 
$$
which computes $R^i {f_2}_* \E$.
 
\subsection{Fibers of $K(\E)$}

   Let $y \in S$ be a closed point, and let $E = \E|_{X_y}$.  Since
$\E$ is flat over $S$ and $R^n {f_2}_* (\E(C)) = 0$ for $n>0$, we have
(cf. ~\cite{mumford.abelian-varieties}, ch.2, n.5, coroll.4,3,2
respectively):
\begin{gather*}
     H^n(X,E(C))=0 \text{ for } n>0 , \\
     (R^0 {f_2}_* \E(C)) \otimes k(y) = H^0(X, E(C)), \\
     R^0  {f_2}_* \E(C) \text{ is locally free.} 
\end{gather*}
and for $F=E(C)|_{C}$ we have
\begin{gather*}
     H^n(C,F(\xi))=0  \text{ for } n>0, \\
     (R^0 {f_1}_* \F(\xi)) \otimes k(y) = H^0(C, F(\xi)), \\
     R^0  {f_1}_* \F(\xi) \text{ is locally free.} 
\end{gather*}
Let $G=F(\xi)|_{\xi}$. It is clear that 
\begin{gather*}
     H^n(\xi,G)=0  \text{ for } n>0, \\
     (R^0 {f_0}_* \G \otimes k(y) = H^0(\xi,G), \\
     R^0  {f_0}_* \G \text{ is locally free.} 
\end{gather*}

It proves that 
$$
   K(\E) \tensor k(y) \simeq K({E_y})
$$

\subsection{Construction of Brill-Noether loci}
If $K$ is a complex 
$$
   E_1 \strelka{\sigma_1}
   E_2 \strelka{\sigma_2}
   E_3
$$
of locally free sheaves on a scheme $S$, we define a closed subscheme
$S_k(K)$ as 
$$
   S_k = \union_{k_1 + k_2 - \rk E_1 = k}
         Z(\Lambda^{k_1} \sigma_1)
         \intersect
         Z(\Lambda^{k_2} \sigma_2)
$$
where $Z(\phi)$ is a scheme of zeroes of $\phi$.

The points $y$ of $S_k$ are exactly the points of $S$
for which $\dim_{k(y)} H(K \tensor k(y)) \ge k$.

Thus, having a flat family of sheaves en a surface (or a scheme) $X$
with base $S$ and using the complex $K(\E)$ constructed in the
previous paragraph we can define the Brill-Noether subschemes $BN^i_k$
such that for every $y \in BN^i_k$ we have $H^i(X,E_y) \ge k$.

If $(X,H)$ is a polarized K3 surface and $\E$ is such a family that
$(c_1(E_y),H)> 0$, then $H^2(X,E_y) = 0$ for every $y \in S$ and
Riemann-Roch theorem gives $BN^1_k = BN^0_{\chi(E_y) - k}$.  In this
situation we define $BN_k= BN^0_k$.

%      no extra [], spell-checked
%      \section{Construction of Brill-Noether loci}
%         \subsection{Complex $K(E)$}
%         \subsection{Relative case}
%         \subsection{Fibers of $K(\E)$}
%         \subsection{Construction of Brill-Noether loci}
%      -------

% ----------------------------------------------------------------------

\section{Construction of the correspondence}
\subsection{The space of Brill-Noether pairs}
\subsubsection{$\alpha$-stable Brill-Noether pairs}
\newcommand{\ep}{\epsilon}
\providecommand{\al}{\alpha}
    Let $E$ be a coherent sheaf on algebraic variety $X$ and $V
\subset H^0(E)$ be a linear subspace. $(E,V)$ is called a
Brill-Noether pair.

Let $H$ be an ample line bundle on $X$ and let $P(E) = P_H(E)$ be the
Hilbert polynomial of $E$ with respect to $H$. Assume that $\rk E \ge
1$.

 {\bf Definition (Le Potier.)}

  (1) Let $\al$ be a polynomial with positive rational
coefficients. Then
\begin{gather*}
 p(E) = \frac{P(E)}{\rk E}, \\
 p_{\al}(E,V)= p(E) + \al \cdot \frac{\dim V}{\rk E} 
\end{gather*}

  (2) The pair $(E,V)$ is called $\al$-stable if $E$ is pure and for
any sub-pair $(F,W) \subset (E,V)$ one has $p_{\al}(F,W) <
p_{\al}(E,V)$.

{\bf Theorem} (Le Potier - He Min): There is a (coarse) moduli space
of $\al$ - stable Brill-Noether pairs denoted as $Syst_{X,\al}(l,
p(E))$, where $l = \dim V$. (cf. \cite{he-min})

\subsubsection{$\epsilon$-stable Brill-Noether pairs}

\begin{thing}
{\bf Lemma.} Let $(E,V)$ be a Brill-Noether pair on a surface $X$.
Assume that $\rk E \ge 2$. Let $\ep$ be a positive rational number
such that $\ep < \frac{1}{\dim V(\rk E - 1)}$. If $(E,V)$ is
$\ep$-stable, then $E$ is semistable.
\end{thing}

{\bf Proof.} Let $(E,V)$ be an $\ep$-stable pair for some positive
rational $\ep$, where $\rk E \ge 2$. Let $F$ be a subsheaf of
$E$. Since $E$ is pure of rank at least 2, we have $\rk F > 0$. Let $e
= \rk E$ and $f = \rk F$. Since $(F,0)$ is a sub-pair of $(E,V)$, we
have $$ p(F) < p(E) + \ep \frac{\dim V}{\rk E}, $$ or \begin{equation}
   \frac{\chi(F_H)}{f} m + \frac{\chi(F)}{f}  < 
   \frac{\chi(E_H)}{e} m + \frac{\chi(E)}{e} + \ep \frac{\dim V}{e}
   \label{eq:(F,0)vs(E,V)}
\end{equation}
It implies that either $\frac{\chi(F_H)}{f} < \frac{\chi(E_H)}{e}$,
in which case $p(F) < p(E)$ and $F$ does not destabilize $E$, or
$\frac{\chi(F_H)}{f} = \frac{\chi(E_H)}{e}$,
in which case 
~\eqref{eq:(F,0)vs(E,V)}
 implies 
$$
  \frac{\chi(F)}{f} - \frac{\chi(E)}{e} < \ep \frac{\dim V}{e},
$$
or
$$
  \chi(F)e - \chi(E)f < \ep f \dim V 
$$ 
In particular, if $\ep f \dim V < 1$, then $\chi(F)e - \chi(E)f \le
0$, which implies $p(F) \le p(E)$. It follows that if $\ep$ is such
that $\ep < \frac{1}{f \dim V}$ for any potentially destabilizing $F$,
then $E$ is semistable. In the definition of stability of a sheaf $E$
it is enough to consider the subsheaves $F \subset E$ for which $\dim
\Supp E/F = \dim \Supp E$, and therefore we may assume that $f \le
\rk(E) - 1$. It implies that if $\ep < \frac{1}{\dim V (\rk E - 1)}$,
then $E$ is semistable.

\begin{thing}
{\bf Lemma.} Let $X$ be a polarized algebraic surface, and let $(E,V)$
be a Brill-Noether pair on $X$. Assume that $\rk E \ge 2$. Let $\ep$
be a positive rational number such that $\ep < \frac{1}{\dim V(\rk E -
1)}$. If $E$ is stable, then the pair $(E,V)$ is $\ep$-stable.
\end{thing}

{\bf Proof.} Assume that $E$ is stable and let $(F,W)$ be a sub-pair of
$(E,V)$.  Let $f = \rk F$ and $e = \rk E$; we have $f > 0$.

Since $E$ is stable, $p(F) < p(E)$, or
\begin{equation}
   \frac{\chi(F_H)}{f} m + \frac{\chi(F)}{f}  < 
   \frac{\chi(E_H)}{e} m + \frac{\chi(E)}{e} 
   \label{eq:FvsE}
\end{equation}

Let $\ep$ be a positive number. The condition that $(F,W)$ does not
$\ep$-destabilize $(E,V)$  is equivalent to the
\begin{equation}
   \frac{\chi(F_H)}{f} m + \frac{\chi(F)}{f} + \ep \frac{\dim W}{f} < 
   \frac{\chi(E_H)}{e} m + \frac{\chi(E)}{e} + \ep \frac{\dim V}{e}
   \label{eq:(F,W)vs(E,V)}
\end{equation}

Now ~\eqref{eq:FvsE} implies that either $\frac{\chi(F_H)}{f} <
\frac{\chi(E_H)}{e}$, in which case ~\eqref{eq:(F,W)vs(E,V)} holds
automatically, or that $\frac{\chi(F_H)}{f} = \frac{\chi(E_H)}{e}$ and
$ \frac{\chi(E)}{e} - \frac{\chi(F)}{f} > 0,$ in which case
~\eqref{eq:(F,W)vs(E,V)} is equivalent to the condition
$$
   \left( f \chi(E) - e \chi(F)\right) + \ep 
   \left( f \dim V - e \dim W \right)    > 0
$$
Since $ f \chi(E) - e \chi(F) \ge 1$, the last condition is equivalent to
the condition
$\ep ( f \dim V - e \dim W ) > -1 $, or 
$$
   \ep < \frac{1}{e \dim W - f \dim V}
$$
In particular, if $\ep < \frac{1}{(e - 1) \dim V}$,
\begin{for-me}
 (which corresponds to the case in which $E$ has a subsheaf $F$ of
rank one such that $V \subset H^0(F)$)
\end{for-me}
then $(E,V)$ is $\ep$-stable.

\begin{for-me}
{\bf} {\bf Note for myself.} Let $v \in K_0(X)$, and let $A(l,v)$ be
the moduli space of $\ep$-stable pairs. It follows that there is a
morphism $A(l,v) \to \bar{M}(v)$, or that there is an open subscheme
$A^{s}(l,v)$ and a morphism $A^{s}(l,v) \to M(v)$.
\end{for-me}

\subsection{The space of good Brill-Noether pairs}
  Let $X$ be an algebraic surface, $r \ge 1$, $v_1 \in \Num(X)$ and
$v_2 \in 1/2 \ZZ$. Let $v=(r,v_1,v_2)$ and let $M=M_H(v)$ be the
moduli space of stable {\it acceptable} sheaves with Chern character
$v$. Let $r'$ be an an integer such that $0 \le r' < r$, and let
$v'=(r',v_1,v_2)$. Let $l=r-r'$, and let $A_l(v) = A_{r,r'}(v)$ be the
open subscheme in $Syst_{\epsilon}(r-r',v)$ which parametrize the
pairs $(E,V)$ such that 
\begin{enumerate}
  \item  $E$ is acceptable, \\ 
  \item  $e_V: V \tensor \OO_X \to E$ is monomorphic, \\
  \item  $E' = \coker e_V$ is acceptable, \\
  \item  $\coker e_V$ is stable.
\end{enumerate}

We call $A_l(v)$ the moduli space of {\it good} Brill-Noether pairs. 
A point $(E,V)$ in the $A_l(v)$ gives an exact sequence
$$
   \triple{V \tensor \OO_X}{E}{E'}
$$

Let $M'=M_H(v')$ be the moduli space of stable {\it acceptable}
sheaves with Chern character $v' = (r',v_1,v_2)$.  The are natural
maps $\pi_1: A_l(v) \to M$ and $\pi_2: A_l(v) \to M'$ which take the
pair $(E,V)$ to $E$ and $E'$ correspondingly:
\begin{equation}
      \correspondence{A_l(v)}{M(v)}{M(v')}{\pi_1}{\pi_2}
      \label{diagram:correspondence}
\end{equation}
\subsection{Description of fibers}

\begin{thing}
\label{description-of-fibers}
{\bf Lemma.}

(1) The fiber of $\pi_1$ at a point $[E] \in M$ is isomorphic to an
open subscheme in the Grassmanian variety $Gr(l,H^0(E))$. 

(2) The fiber of $\pi_2$ at a point $[E'] \in M'$ is isomorphic to an
open subscheme in $Gr(l,H^1(E')\dual)$.
\end{thing}

The proof of (1) is clear: the fiber of $\pi_1$ is an open subscheme
in $Gr(l,H^0(E))$ parametrizing linear subspaces $V$ for which the
evaluation map $e_V$ is monomorphic and $\coker e_V$ is acceptable and
stable.

\label{fibers-of-pi_2}

To describe fibers of $\pi_2$ at $[E'] \in M'$, we proceed in a few
steps.

Step 1: fix a vector space $W$ of dimension $r-r'$ and define the
moduli space $A'_{r,r'}$ of ``rigified'' Brill-Noether pairs which
parametrizes data $(E,i:W \to H^0(E))$, where $i$ is monomorphic and
$\dim W = r-r'$. In other words, $A'_{r,r'}$ is a moduli space of
triples $(E, V \subset H^0(E), j: V \isom W)$. There is a forgetting
morphism $p: A'_{r,r'} \to A_{r,r'}$, and it is clear that $A'_{r,r'}$
is a principle $GL(W)$-bundle over $A_{r,r'}$.

Let $\pi'_2: A' \to M'$ take $(E, W)$ to $E':= \coker i_W$. There
is a commutative diagram 
$$ 
   \smallXyMatrix { {A'_{r,r'}} \ar[d]_{\alpha} \ar[rdd]^{\pi'_2}& {}
   \\ {A_{r,r'}} \ar[rd]_{\pi_2} & {} \\ {} & {M'} }
$$

  A point $(E,i)$ in $A'$ gives an extension
\begin{equation}
  \triple{W \tensor \OO_X}{E}{E'}
  \label{extension-by-W}
\end{equation}
where $E' = [\coker i_W]$.

\begin{for-me}
------------------------------------------------
{\bf (2)}
Zdes' nachinaetsya {\bf oshibka}: ya zabyvayu, chto $M'$ parametrizuet
klassy izomorfizma.:

It is clear that the fiber of $\pi'_2$ over 
$[E'] \in M'$ is an open subscheme in $\Ext^1(E',\triv{W})$
given by the condition that the middle term $E$ in the extension
~\eqref{extension-by-W} is acceptable and stable. 

  The fiber $(\pi'_2)^{-1}(E')$ is invariant under the $GL(W)$ -action
and therefore
$(\pi_2)^{-1}(E') \simeq (\pi'_2)^{-1}(E') / GL(W)$ 

[Zdes', vidimo, proishodit sleduyuschee: Esli $p: X \to Y$ - $G$ -
rassloenie i $Z$ - zamknutaya $G$ - invariantnaya podshema v $X$, to
$p(Z)$ zamknuto v $Y$ (dazhe hotya morfizm ne proektiven), i
predstavlyaet funktor orbit].

(konec oshibochnogo rassuzhdeniya)
------------------------------------------------------------
\end{for-me}

% (we assume that $h^0(X,\ox) = 1$.)

Let $m' = [E'] \in M'$, and let $F_{m'}$ be the fiber of
$\pi_2$ over $m'$. 
% Let us fix a vector space $W$ of dimension $r-r'$.
Let $\Ext^1(E', W \tensor \ox)^{a}$ be the open subscheme in
$\Ext^1(E', W \tensor \ox)$ classifying stable acceptable extensions
There is a natural morphism
$$
   \theta: \Ext^1(E', W \tensor \ox)^{a} \to F_{m'}
$$
which takes  the extension class
$$
   \xyTriple{W \tensor \ox}{E}{E'}{\alpha}{\beta}
$$ to the pair $(E,h^0(\alpha)(W))$, where $h^0(\alpha): W \to
H^0(X,E)$ is induced by $\alpha$.

It is clear that $\theta$ is epimorphic.
\begin{for-me}
Indeed, if $(E,V) \in F_{m'}$, then there is an exact sequence
$$
   \xyTriple{V \tensor \ox}{E}{\coker e_V}{e_V}{can}
$$
and an isomorphism $i: \coker e_V \isom E'$. The induced extension of
$E'$ gives a class $e \in \Ext^1(E', W \tensor \ox)$.
\end{for-me}

Step 2: Note that there is an action of the group $GL(W) \times
\Aut(E')$ on $\Ext^1(E',W \tensor \ox)$ which preserves the (subscheme
of) stable acceptable extensions: if $(\gamma, \delta) \in GL(W)
\times \Aut(E')$, then we define $(\gamma,\delta) \circ e$ to be the
bottom row of the diagram
$$
  \xymatrix
  {
    0 \ar[r] & {W \tensor \ox} \ar[r]^-{\alpha}\xyIsom[d]_{\gamma} & E \ar[r]^-{\beta}\xyRavno[d] & E' \ar[r]\xyIsom[d]_{\delta} & 0 \\
    0 \ar[r] & {W \tensor \ox} \ar[r]^-{\alpha\gamma^{-1}} & E \ar[r]^-{\delta\beta} & E' \ar[r] & 0 \\
  }
$$

It is clear that $\theta$ is invariant under this action,
i.e., $\theta( (\gamma,\delta) \circ e ) = \theta(e)$.

\begin{thing}Lemma.
If $e_i \in \Ext^1(W \tensor \ox,E')^a$, $i = 1,2$, and $\theta(e_1) =
\theta(e_2)$, then $e_2 = (\gamma,\delta) \circ e_1$ for some element
$(\gamma,\delta) \in GL(W) \times \Aut(E')$.
\end{thing}

Indeed, let $V_i = h^0(\alpha_i)(W)$, $i = 1,2$; since $\theta(e_1) =
\theta(e_2)$, there is an isomorphism $i: E_1 \to E_2$ which takes
$V_1$ to $V_2$. Let us construct $\gamma \in GL(W)$ from the
diagram
$$
   \xymatrix
   {
      W \ar[r]^{\alpha_1}\ar@{-->}[d]^{\gamma} & V_1\ar[d]^{i} \\
      W \ar[r]^{\alpha_2} & V_2 
   }
$$

Since $i: E_1 \to E_2$ takes $\Im \alpha_1$ to $\Im \alpha_2$, the
isomorphism $i$ induces an isomorphism $\delta: E' \to E'$ such that
the following diagram is commutative:
$$
  \xymatrix
  {
    0 \ar[r] & {W \tensor \ox} \ar[r]^-{\alpha_1}\ar[d]^-{\gamma} & E_1 \ar[r]^-{\beta_1}\ar[d]^-{i} & E' \ar[r]\ar[d]^-{\delta} & 0 \\
    0 \ar[r] & {W \tensor \ox} \ar[r]^-{\alpha_2} & E_2 \ar[r]^-{\beta_2} & E' \ar[r] & 0 \\
  }
$$
But now, twisting with $(\gamma,\delta)$, we get an equivalence of
extensions $(\gamma,\delta) \circ e_1 \equiv e_2$:
$$
  \xymatrix
  {
    0 \ar[r] & {W \tensor \ox} \ar[r]^-{\alpha_1\gamma^{-1}}\xyRavno[d] & E_1 \ar[r]^-{\delta\beta_1}\xyIsom[d]_-{i} & E' \ar[r]\xyRavno[d] & 0 \\
    0 \ar[r] & {W \tensor \ox} \ar[r]^-{\alpha_2} & E_2 \ar[r]^-{\beta_2} & E' \ar[r] & 0 \\
  }
$$

Step 3: Since $E'$ is stable, it is also simple and the action of the
group $GL(W) \times \Aut(E')$ reduces to the action of $GL(W)$.  Let
$\alpha_e$ be the image of $e$ under the isomorphism $u:
\Ext^1(E',\triv{W}) \to \Hom(W\dual, \Ext^1(E',\OO))$. Lemma
~\ref{stability-of-extensions-by-W-tensor-O} implies that if
$\alpha_e$ is not monomorphic, then $E$ is not stable . Consider the
space of monomorphic linear maps $\Monom(W\dual,\Ext^1(E',\OO))$. It
has free $GL(W)$-action, and $\Monom(W\dual,\Ext^1(E',\OO)) / GL(W)
\simeq \Gr(l,\Ext^1(E',\OO))$, where $l = \dim W = r-r'$. It follows
that the fiber of $\pi_2$ is an open subscheme in the Grassmanian
variety $\Gr(r-r',\Ext^1(E',\OO))$ given by the condition that $E$ in
the extension defined by $\alpha_e \in \Monom(W\dual,\Ext^1(E',\OO))$
is acceptable and stable.

\subsection{Nonemptiness of fibers}
\begin{thing}
\label{fiber-of-pi1-in-gg-PicZ-case}
{\bf Lemma.} Assume that $[E] \in M$ and either $\rk E = 1$, or $\rk E
\ge 2$ and the following two conditions are true:

(a) $E$ is globally generated, and

(b) $\Pic X \simeq \ZZ, c_1(E) = 1$.

Then the fiber of $\pi_1$ at $[E] \in M$ is a {\it dense} open
subscheme in the Grassmanian variety $Gr(l,H^0(E))$.
\end{thing}

{\bf Proof.} In case $\rk E = 1$ the lemma is easy. If $\rk E \ge 2$,
then, first, the condition (a) and Lemma ~\ref{injectivity-lemma}
imply that the generic evaluation map is monomorphic. Second, the same
condition and Lemma ~\ref{good-defeneration-in-the-gg-case} imply
that the generic factorsheaf of the evaluation map is
acceptable. Finally, condition (b) and the Lemma
~\ref{stability-of-WO-factors} imply that the generic factorsheaf is
stable.

\begin{thing}
\label{fiber-of-pi2-in-gg-PicZ-case}
{\bf Lemma.} Assume that $[E'] \in M'$ and the following two
conditions are true:

\begin{itemize}
  \item[(a)] $\Pic X \simeq \ZZ, c_1(E) = 1$,

  \item[(b)] Either one of the following is true:
      \begin{itemize}
          \item[(1)] $\rk F \ge 2$;
          \item[(2)] $\rk F = 1$, $F \simeq J_{\xi}(L)$, where $L$ is ample,
$\xi$ is simple and the pair $(\xi, L+K_X)$ is Caley-Bacharash;
          \item[(3)] $\rk F = 0$, $F \simeq (i_C)_* B$, where $i_C: C
\to X$ is a reduced irreducible curve on $X$, and $B$ is an invertible
sheaf on $C$.
      \end{itemize}
\end{itemize}

Then the fiber of $\pi_2$ at $[E'] \in M'$ is a {\it dense} open
subscheme in the Grassmanian variety $Gr(l,H^1(E')\dual)$.
\end{thing}
{\bf Proof. } The Lemma
~\ref{criteria-for-a-generic-extension-to-be-acceptable}
%valid ref
and condition (b) imply that the generic extension of $E'$ with $W
\tensor \OO_X$ is acceptable, and condition (a) and Lemma
~\ref{stability-of-extensions-by-W-tensor-O} imply that the generic
extension is stable.

\begin{for-me}
  \input{for-me.abel-jacobi-for-curves}
  % this section is now in the first part of the paper
  % besides, there is a mistake in the computation of Ext^1 
\end{for-me}

%  no extra [], spell-checked
% \section{Construction of the correspondence}
%    \subsection{(*) The space of Brill-Noether pairs}
%       \subsubsection{(*) $\alpha$-stable Brill-Noether pairs}
%       \subsubsection{(..) $\epsilon$-stable Brill-Noether pairs}
%    \subsection{(*) The space of good Brill-Noether pairs}
%    \subsection{Description of fibers}
%    \subsection{Nonemptiness of fibers}
% -----------

%  \input{ composition }
%             composition of correspondences

\section{The correspondence and the Brill-Noether loci}
Let $(X,H)$ be a polarized surface, and assume that $h^1(X,\ox) = 0$.
Let $L = \ZZ \oplus \Num X \oplus \ZZ \cdot \frac{t}{2}$, where $t =
2$ if $\Num X$ is even and $t=1$ otherwise.  There is a Chern
character map $\ch: K_0(X) \to L$, $[E] \mapsto \ch(E) =
(\rk(E),\ch_1(E),\ch_2(E))$.

Let $v \in L$, $v = (r,v_1,v_2)$, and let $(E,V) \in A_l(v)$.
%
%Consider the correspondence ~\eqref{diagram:correspondence} 
%
We get an exact sequence 
\begin{equation}
   \triple{V \tensor \ox}{E}{E'}
   \label{diagram:VO-E-E'}
\end{equation}
and the long exact sequence
$$
   \triple{V}{H^0(X,E)}{H^0(X,E')}
$$ which implies $h^0(X,E') = h^0(X,E) - l$. Let $A_l^k(v)$ be the
preimage of $BN_k(v)$ under $\pi_1$, and let $v'=(r-l,v_1,v_2)$. Since
$c_t(E') = c_t(E)$, the diagram ~\eqref{diagram:correspondence}
induces the diagram
\begin{equation}
   \smallXyMatrix
   {
        && A_l^k(v) 
        \ar[lldd]_{(\pi_1)_k}\ar[rrdd]^{(\pi_2)_k}\xyArrow{\cap}{d} && \\ 
        && A_l(v) \ar[ld]_{\pi_1}\ar[rd]^{\pi_2} && \\ 
        BN_k(v) \xyArrow{\subset}{r} & M(v) && M(v') & 
        BN_{k-l}(v') \xyArrow{\supset}{l}  
   }
   \label{diagram:bn-subcorrespondence}
\end{equation}

It is clear that $((\pi_1)_k)^{-1} (m) = (\pi_1)^{-1} (m)$ and
$((\pi_2)_k)^{-1} (m') = (\pi_2)^{-1} (m')$.

%  no extra [], spell-checked

%  \input{ numerical-structure }
%  ``numerical structure'' is now contained in Part 1

\section{Globally generated line bundles on curves}

\subsection{Varieties ${}_pA^k_d$ and ${}_pB^k_d$}

\subsubsection{}
% {\bf 1.}
 Let $C$ be a smooth irreducible genus $g$ curve and $L$ be a degree
$d$ line bundle on $C$. Let $L^D = K_C L^{-1}$, and let us fix a point
$p \in C$. The adjunction sequence
$$
   \triple{L(-p)}{L}{L_p}
$$ 
induces the long exact sequence
$$
   0 \to H^0(C,L(-p)) \to H^0(C,L) \strelka{e(L,p)} L_p \to 
         H^1(C,L(-p)) \to H^1(C,L) \to 0
$$ Let $M = L(-p)$.

\begin{for-me}
 The sequence above can be rewritten as
$$
   0 \to H^0(C,M) \to H^0(C,M(p)) \to M(p)_p \to 
         H^1(C,M) \to H^1(C,M(p)) \to 0
$$
\end{for-me}

The long exact sequence above implies that there are two possibilities:
either

  (1) $e(L,p)$ is epimorphic, in which case the isomorphism $H^1(C,M)
\isom H^1(C,M(p))$ together with Serre duality imply $e(M^D,p) = 0$;

or

  (2) $e(L,p) = 0$, in which case $e(M^D,p)$ is epimorphic.

   In other words, $L$ is globally generated at $p$ iff $M^D$ is not
globally generated at $p$; and $L$ is not globally generated at $p$
iff $h^0(C,M)=h^0(C,M(p))$.

\subsubsection{}
% {\bf 2.} 
Let now $V_d^k$ be a scheme parametrizing invertible sheaves of degree
$d$ such that $\dim H^0(C,L) = k$; this is a difference of two
Brill-Noether loci in the Picard variety of $C$; and let ${}_pA_d^k$
be the (closed) subscheme in $V_d^k$ parametrizing invertible sheaves
$L$ not globally generated at $p$ (it is constructed in paragraph
~\ref{bundles-not-globally-generated-at-a-given-point}.) We define
${}_pB_d^k$ as ${}_pB_d^k = D({}_pA_{2g-2-d}^{k-d-1+g})$; in other
words, $L \in B$ iff $L^D \in A$.

\begin{for-me}
(Is there a global locally closed subscheme ${}_pA_d$ in 
$\Pic^d C$ giving ${}_pA_d^k$ on each $V_d^k$?)
\end{for-me}

The map $L \mapsto L(p)$ establishes isomorphisms
\begin{gather}
    V_d^k - {}_pB_d^k \isom {}_pA_{d+1}^k, \\
    {}_pB_d^{k-1} \isom V_{d+1}^k -  {}_pA_{d+1}^k
\end{gather}

In other words, the first isomorphism can be formulated as follows:
twisting with $\oc(p)$ establishes an isomorphism of the (open)
subscheme in moduli of $M$ given by the condition $h^0(C,M(p)) =
h^0(C,M)$ and the (closed) subscheme in the moduli of $L$ given by the
condition that $L$ is not globally generated at $p$.

The isomorphisms of twisting with $p$ and $-p$ are shown on the
Figures ~\ref{figure:twisting-with-p} and
~\ref{figure:twisting-with-p} as acting in the plane with
$(d,k)$-coordinates.

\begin{figure}[!h]   % here, literally
\begin{center}
   \setlength{\unitlength}{3947sp}%
\begingroup\makeatletter\ifx\SetFigFont\undefined%
\gdef\SetFigFont#1#2#3#4#5{%
  \reset@font\fontsize{#1}{#2pt}%
  \fontfamily{#3}\fontseries{#4}\fontshape{#5}%
  \selectfont}%
\fi\endgroup%
\begin{picture}(3459,2467)(1039,-4820)
\thinlines
{\color[rgb]{0,0,0}\put(1951,-3511){\vector( 2,-1){2400}}
}%
{\color[rgb]{0,0,0}\put(1951,-3211){\vector( 1, 0){2400}}
}%
\put(1179,-3226){\makebox(0,0)[lb]{\smash{\SetFigFont{5}{6.0}{\familydefault}{\mddefault}{\updefault}{\color[rgb]{0,0,0}$V_d^k - {}_pB_d^k$}%
}}}
\put(1449,-3526){\makebox(0,0)[lb]{\smash{\SetFigFont{5}{6.0}{\familydefault}{\mddefault}{\updefault}{\color[rgb]{0,0,0}${}_pB_d^k$}%
}}}
\put(4498,-3243){\makebox(0,0)[lb]{\smash{\SetFigFont{5}{6.0}{\familydefault}{\mddefault}{\updefault}{\color[rgb]{0,0,0}${}_pA_{d+1}^k$}%
}}}
\put(4486,-3526){\makebox(0,0)[lb]{\smash{\SetFigFont{5}{6.0}{\familydefault}{\mddefault}{\updefault}{\color[rgb]{0,0,0}$V_{d+1}^k - {}_pA_{d+1}^k$}%
}}}
\put(4478,-4471){\makebox(0,0)[lb]{\smash{\SetFigFont{5}{6.0}{\familydefault}{\mddefault}{\updefault}{\color[rgb]{0,0,0}${}_pA_{d+1}^{k+1}$}%
}}}
\put(4463,-4771){\makebox(0,0)[lb]{\smash{\SetFigFont{5}{6.0}{\familydefault}{\mddefault}{\updefault}{\color[rgb]{0,0,0}$V_{d+1}^{k+1} - {}_pA_{d+1}^{k+1}$}%
}}}
\put(1039,-2547){\makebox(0,0)[lb]{\smash{\SetFigFont{12}{14.4}{\familydefault}{\mddefault}{\updefault}{\color[rgb]{0,0,0}${\rm Pic}^d \; C$}%
}}}
\put(4187,-2521){\makebox(0,0)[lb]{\smash{\SetFigFont{12}{14.4}{\familydefault}{\mddefault}{\updefault}{\color[rgb]{0,0,0}${\rm Pic}^{d+1} \; C$}%
}}}
\end{picture} 
   \caption{Twisting with $p$}
   \label{figure:twisting-with-p}
\end{center}
\end{figure}
\vspace*{3mm}
\begin{figure}[!h]   % here, literally
\begin{center}
    \setlength{\unitlength}{3947sp}%
\begingroup\makeatletter\ifx\SetFigFont\undefined%
\gdef\SetFigFont#1#2#3#4#5{%
  \reset@font\fontsize{#1}{#2pt}%
  \fontfamily{#3}\fontseries{#4}\fontshape{#5}%
  \selectfont}%
\fi\endgroup%
\begin{picture}(3525,2167)(976,-4520)
\thinlines
{\color[rgb]{0,0,0}\put(4351,-4411){\vector(-1, 0){2400}}
}%
{\color[rgb]{0,0,0}\put(4351,-4111){\vector(-2, 1){2400}}
}%
\put(4478,-4471){\makebox(0,0)[lb]{\smash{\SetFigFont{5}{6.0}{\familydefault}{\mddefault}{\updefault}{\color[rgb]{0,0,0}${}_pA_{d+1}^{k+1}$}%
}}}
\put(1039,-2547){\makebox(0,0)[lb]{\smash{\SetFigFont{12}{14.4}{\familydefault}{\mddefault}{\updefault}{\color[rgb]{0,0,0}${\rm Pic}^d \; C$}%
}}}
\put(4187,-2521){\makebox(0,0)[lb]{\smash{\SetFigFont{12}{14.4}{\familydefault}{\mddefault}{\updefault}{\color[rgb]{0,0,0}${\rm Pic}^{d+1} \; C$}%
}}}
\put(4501,-4261){\makebox(0,0)[lb]{\smash{\SetFigFont{5}{6.0}{\familydefault}{\mddefault}{\updefault}{\color[rgb]{0,0,0}$V_{d+1}^{k+1} - {}_pA_{d+1}^{k+1}$}%
}}}
\put(4501,-3361){\makebox(0,0)[lb]{\smash{\SetFigFont{5}{6.0}{\familydefault}{\mddefault}{\updefault}{\color[rgb]{0,0,0}$V_{d+1}^k - {}_pA_{d+1}^k$}%
}}}
\put(4501,-3511){\makebox(0,0)[lb]{\smash{\SetFigFont{5}{6.0}{\familydefault}{\mddefault}{\updefault}{\color[rgb]{0,0,0}${}_pA_{d+1}^k$}%
}}}
\put(1276,-4186){\makebox(0,0)[lb]{\smash{\SetFigFont{5}{6.0}{\familydefault}{\mddefault}{\updefault}{\color[rgb]{0,0,0}${}_pB_d^{k+1}$}%
}}}
\put(976,-4411){\makebox(0,0)[lb]{\smash{\SetFigFont{5}{6.0}{\familydefault}{\mddefault}{\updefault}{\color[rgb]{0,0,0}$V_d^{k+1} - {}_pB_d^{k+1}$}%
}}}
\put(1351,-2911){\makebox(0,0)[lb]{\smash{\SetFigFont{5}{6.0}{\familydefault}{\mddefault}{\updefault}{\color[rgb]{0,0,0}${}_pB_d^{k-1}$}%
}}}
\put(1276,-3211){\makebox(0,0)[lb]{\smash{\SetFigFont{5}{6.0}{\familydefault}{\mddefault}{\updefault}{\color[rgb]{0,0,0}$V_d^k - {}_pB_d^k$}%
}}}
\end{picture} 
    \caption{Twisting with $-p$}
    \label{figure:twisting-with-minus-p}
\end{center}
\end{figure}

% {\bf 3.}
\subsubsection{}
\newcommand{\uA}{\underline{{A}}}
\newcommand{\uB}{\underline{{B}}}
We will now vary the point $p$. In paragraph
~\ref{bundles-not-globally-generated-at-some-point} we construct a 
scheme $\uA_d^{k}$ (and $\uB_d^{k}$) such that there are diagrams
$$
   \smallXyMatrix
   {
      V_d^k \times C \ar[d] \xyArrow{\supset}{r}  & {\uB_d^k} \ar[ld]  \\
      C  &  {} 
   }
   \, \, \, \, \, \, \, \, \, \, \, \, 
   \smallXyMatrix
   {
      V_d^k \times C \ar[d] \xyArrow{\supset}{r}  & {\uA_d^k} \ar[ld]  \\
      C  &  {} 
   }
$$ such that the fibers are isomorphic to ${}_pB_d^k$ and ${}_pA_d^k$,
respectively.

  Let $B_d^k \subset V_d^k$ be the image of $\uA_d^k$ under the of
projection map $\uB_d^k \to V_d^k$. It can be described as a closed
subscheme in $V_d^k$ parametrizing such $M$ that $h^0(M)= h^0(M(p))$
for some point $p \in C$. In the same way there is a closed subscheme
$A_d^k \subset V_d^k$ parametrizing non-globally generated invertible
sheaves.

\subsubsection{}
%{\bf 4.} 
Now consider the morphism
$$
    \phi: \;
    (V_{d-1}^k \times C)  - \uB_{d-1}^k 
    \; \; \to \; \; 
    A_d^k
$$ mapping $(M,p)$ to $M(p)$. It is clear that this morphism is
surjective and for $k \ge 1$ is quasi-finite (the set of points where
the given line bundle is not globally generated is finite).

The map $\phi$ can be included into the diagram
$$
   \smallXyMatrix
   {
      (V_{d-1}^k \times C) - \uB_{d-1}^k   \ar[rd]_{\phi}\xyIsom[r] &
      {\uA_{d}^k}                          \ar[d]^{\pi}            \\
      {}  &  A_d^k 
   }
$$
where $\pi$ is quasi-finite.

{\bf Corollary.} Let $k \ge 1$. Then for each irreducible component
$(A_d^k)_i$ of $A_d^k$ (if any) we have $\dim (A_d^k)_i \le \max_j
\dim (V_{d-1}^k)_j + 1$, where $(V_{d-1}^k)_j$ are the irreducible
components of $V_{d-1}^k$.

% what are the theorems on the connectedness / irreducibility / 
% equidimensionality ?
% $V$ can be reducible even if $\dim V > 0$: example - $g = 5$.

% is it finite?
% fiber at L = {points p at which L is not gg = H^0(C,L) vanishes;
% generically speaking, can vary with L. 
% but this is fine; consider the normalization of a curve with nodes, f.e.

\subsection{Varieties of non-globally generated line bundles.}

\subsubsection{Bundles not globally generated at a given point}
\label{bundles-not-globally-generated-at-a-given-point}
   Let $C$ be a curve and $L$ be a family of invertible sheaves on $C$
with base $T$, i.e., an invertible sheaf $L$ on $C \times T$. 

   In this paragraph we construct a subscheme of bundles not globally
generated at a given point $p \in C$.

%(We will be most interested in the universal case in which $T = \Pic^d C$ and
% $L$ is the universal line bundle).

 Let $\pi: C \times T \to T$ be the projection. Let us fix a point $p
\in C$, and let $D_p = p \times T$. Applying $R^1 \pi_*$ to the
adjunction sequence on $C \times T$
$$
   0 \to L(-D_p) \to L \to L|_{D_p} \to 0
$$
we get the exact sequence
$$
   R^1 \pi_* L(-D_p) \strelka{\phi_p} R^1 \pi_* L \to 0
$$
(we consider the $R^1$ piece since base change works better for top
cohomology groups.)

Let $t\in T$ and $L_t = L|_{C \times t}$. Tensoring the exact sequence
above with the residue field $k(t)$, we get
$$
   H^1(C,L_t(-p)) \strelka{(\phi_p)(t)} H^1(C,L_t) \to 0   
$$    
which is a part of the long cohomological sequence associated with
$$
   0 \to (L_t)(-p) \to L_t \to (L_t)|_p \to 0
$$

In particular, $L_t$ is not globally generated at $p \in C$ iff
$\ker (\phi_p)(t) \ne 0$.

In general, for a morphism of coherent sheaves $\phi: F \to G$ the set
of points $t \in T$ such that $\phi(t)$ is not an isomorphism is
constructible.

\begin{for-me}
 for $T$ can be decomposed into a union of locally closed
subschemes $T_{i}$ such that on every strata $T_{i}$ the sheaves $F$
and $G$ are locally free. The degeneration set of $\phi|_{T_{i}}$ is
closed in $T_{i}$ and therefore is locally closed in $T$, and it
follows that the set of points $t \in T$ for which $\phi \tensor k(t)$
is not an isomorphism is constructible.

   It follows that there is a constructible subset $A_pT$ in $T$
parametrizing points $t \in T$ such that the invertible sheaf $L_t$ on $C$
is not globally generated at $p$.
\end{for-me}

Now let $V^kT$ be a (locally closed) subscheme in $T$ given by the
condition $\dim_k H^0(C,L_t) = k$. The sheaf $R^1 \pi_* L$ is locally
free when restricted to $V^kT$, and tensoring the exact sequence
$$
    0 \to \ker \phi_p \to R^1 \pi_* L(-D_p) \strelka{\phi_p} R^1 \pi_*
    L \to 0
$$ with $k(t)$ demonstrates that the set of $t \in V^k T$ such that
$L_t$ is not globally generated at $p \in C$ is the support of the
coherent sheaf $\ker \phi_p$ which implies that it is closed.

\subsubsection{Varying the point $p \in C$.}
\label{bundles-not-globally-generated-at-some-point}
Let $C$ be a curve and $L$ be a family of invertible sheaves on $C$
with a base $T$.  In this paragraph we construct a (relative) variety
of bundles not globally generated at a some point $p \in C$.

 Let $C_1 = C_2 = C$, and consider the diagram of projection morphisms
$$
   \xymatrix
   {
      C_1 \times C_2 \times T \ar[d]^{\pi} \\
      C_1 \times T \ar[d]^{\pi'} \\
      T
   }
$$
Let $\Delta$ be a divisor on $C_1 \times C_2 \times T$ obtained by
pulling back the diagonal from $C_1 \times C_2$, and abusing notations we
write $L$ instead of $pr_{C_2 \times T}^* L$.

Applying $R^1\pi_*$ to the adjunction sequence
$$
   0 \to L(-\Delta) \to L \to L|_{\Delta} \to 0,
$$
we get the exact sequence
$$
   R^1 \pi_* L(-\Delta) \strelka{\phi} R^1 \pi_* L \to 0
$$
of coherent sheaves on $C_1 \times T$. If $(p,t) \in C_1 \times T$,
then, tensoring the exact sequence above with the residue filed
$k(p,t)$ of the point $(p,t) \in C_1 \times T$, we get the exact
sequence
$$
   H^1(C,L_t(-p)) \strelka{\phi \tensor k(p,t)} H^1(C,L_t) \to 0   
$$ as in the previous paragraph, and it follows that the constructible
subset $\Deg \phi \subset C_1 \times T$ consisting of points $(p,t)$
for which $\phi \tensor k(p,t)$ is not a fibervise isomorphism
coincides with the set of points $(p,t) \in C_1 \times T$ for which
$L_t$ is not globally generated at $p \in C$. (Now taking the direct
image $\pi'(\Deg \phi)$, we can get a subset $A(T)$ in $T$ of points
$t \in T$ such that the invertible sheaf $L_t$ is not globally
generated.)

\begin{for-me}
% I never use this
\subsubsection{Relative case.}
Let $S$ be a scheme, $C/S$ be a family of curves over $S$, and $L$ be
a relative family of invertible sheaves on $C/S$ with the base $T/S$,
i.e., an invertible sheaf over $C \fibered_S T$. (We will be most
interested in the universal case when $T/S = \Pic^d (C/S)$ is a
relative Picard scheme.  Let $C_1 = C_2 = C$, and consider the diagram
of projections
$$
   \xymatrix
   {
      C_1 \fibered_S C_2 \fibered_S T \ar[d]^{\pi} \\
      C_1 \fibered_S T \ar[d]^{\pi'} \\
      T \ar[d] \\
      S
   }
$$
Let $\Delta$ be the preimage of the diagonal in $C_1 \fibered_S C_2$,
and $L$ be the preimage of the family $L$ on $C_2 \fibered_S T$.

Applying $\pi_*$ to the adjunction sequence
$$
   0 \to L(-\Delta) \to L \to L|_{\Delta} \to 0
$$
we get the exact sequence
$$
   R^1 \pi_* L(-\Delta) \strelka{\phi} R^1 \pi_* L \to 0
$$
of coherent sheaves on $C_1 \fibered_S T$. If $s \in S$ and $(p,t) \in
C_1 \fibered_S T$ lie over $s \in S$, then, tensoring the exact sequence
above with the residue filed $k(p,t)$ of the point $(p,t) \in C_1
\fibered_S T$, we get the exact sequence
$$
   H^1(C_s,L_t(-p)) \strelka{\phi_s \tensor k(p,t)} H^1(C_s,L_t) \to 0   
$$    
It follows that the set $\Deg \phi$ of points $(p,t)$ such that
$\phi_s \tensor k(p,t)$, $s= \pi''((p,t))$ is not a fibervise
isomorphism is equal to the set of pairs $(p,t)$ such that $L_t$ is
not globally generated at $p \in C_s$.  Note that $\Deg \phi$ is a
constructible subset in $C \fibered_S T$, and $\pi(\Deg \phi)$ is a
constructible subset in $T$ parametrizing non-globally generated
invertible sheaves in the family $T/S$.
\end{for-me}

\subsection{Globally generated line bundles
                 on Brill-Noether-general curves}

{\bf Lemma.} Let $C$ be a curve such that every $V^i_d$ is irreducible
and of the expected dimension $\rho = g- i(i-\chi)$, where $\chi =
d+1-g$.

Assume that we are given $k \ge 2$ such that $\rho = g-k(k-\chi) \ge
0$.  Then there is a globally generated line bundle $L \in V^k_d$.

{\bf Proof.} We prove that $A_d^k$ is equidimensional in $V_d^k$ and
compute its codimension.
  
   Consider the variety $\uA_d^k$ (which might be empty). It has a map
to $C$ with fibers ${}_pA^k_d$. We have ${}_pA^k_d \simeq V^k_{d-1} -
{}_pB^k_{d-1}$, where ${}_pB^k_{d-1}$ is closed in $V^k_{d-1}$. It
follows that ${}_pA^k_d$ is either irreducible of dimension $\dim
V^k_{d-1} = g - k(k-(\chi-1)) = \rho - k$ or empty, which implies that
$\uA_d^k$ is either equidimensional of dimension $k(k-(\chi-1)) + 1$
or empty. (Remark: the cases $\rho < k$, $\rho = k$ and $\rho > k$
should be considered separately.)

    Now there is a quasi-finite map $\uA_d^k \to A_d^k$.  It follows
that $A_d^k$ is equidimensional of dimension $\rho - k + 1$ or empty.

    Now $A_d^k$ is a subscheme in $V_d^k$ which is irreducible of
dimension $\rho$. It follows that $A_d^k$ is of codimension $k-1$ in
$V_d^k$ or empty. In particular, there is a globally generated line
bundle $L \in V_d^k$.

The proof is illustrated on the following diagram:

$$
   \xymatrix
   {
       { {}_pA^k_d }  \xyArrow{\simeq}{r} \ar[d] &  { V^k_{d-1} - {}_pB^k_{d-1}} \\
       { \uA_d^k }  \ar[r] \ar[d] &  C \\
       A_d^k  \xyArrow{\subset}{r} &  V_d^k 
   }
$$

\begin{thing}
{\bf Proposition.} Let $C$ be a generic curve of degree $g$, and
assume that we are given $d$, $k \ge 2$ such that $\rho = g-k(k-\chi)
\ge 0$.  Then there is a globally generated line bundle $L \in V^k_d$.
\end{thing}

The proof follows from the Lazarsfeld's theorem on the dimension of
$W^r_d$ of a generic genus $g$ curve and the previous lemma.

\begin{thing}
{\bf Existence of globally generated vector bundles in moduli spaces
on K3} Let $(X,H)$ be a polarized K3 surface. Assume that $\Pic X
\simeq \ZZ h$, $h = [H]$, and consider a nonempty moduli space $M$ of
vector bundles on $S$ with $\rk E = r$, $c_1(E) = h$ and $c_2(E) =
d$. Then there is a globally generated vector bundle $E \in M$.
\end{thing}

{\bf Proof.} First, by previous Lemma, there is a curve $C$ on $X$ in
the linear system $|h|$ and a globally generated line bundle $L$ on
$C$ of degree $2g-2-d$. Then the Lazarsfeld's construction (cf. the
first part of the paper) implies that there is a globally generated
vector bundle $E$ in the moduli space $M$.

% \no []

% -----------------------------------------------------------------

\end{document}